\documentclass[a4paper,11pt]{article}
\usepackage{amsmath}
\usepackage{amssymb}
\usepackage{texdraw}
\usepackage{mathrsfs}
\usepackage{epsfig}
\usepackage{amsfonts}
\usepackage{graphicx} 
\usepackage{float} 
\usepackage{subfigure} 
\usepackage[dvips]{psfrag}
\usepackage{graphicx}
\usepackage{epstopdf}
\usepackage{subfigure}
\usepackage{color}
\usepackage{booktabs}

\usepackage{geometry}
\usepackage{hyperref}
\hypersetup{colorlinks=true,
            linkcolor=blue,
            anchorcolor=blue,
            citecolor=blue}

\newtheorem{Proposition}{Proposition}[section]

\newtheorem{Assumption}{Assumption}[section]
\newtheorem{Lemma}{Lemma}[section]
\newtheorem{Theorem}{Theorem}[section]

\newtheorem{Remark}{Remark}[section]

\newcommand{\bpf}{{\bf Proof:\ \ }}
\newcommand{\epf}{\mbox{}\hfill $\Box$}

\usepackage{caption}
\usepackage {indentfirst}

\begin{document}
\setlength{\baselineskip}{15.2pt} \pagestyle{myheadings}

\title{\bf The discontinuous planar piecewise linear system with two nodes has at most two limit cycles}

\author{Lu Chen{$^{*}$}, Changjian Liu\\
{\small School of Mathematics (Zhuhai), Sun Yat-sen University,}\\
{ \small  Zhuhai, 519082,  P.R. China }}

\renewcommand{\thefootnote}{}
\makeatletter\def\Hy@Warning#1{}\makeatother
\footnotetext{$^*$Corresponding author.\par E-mail addresses: chenlu73@mail2.sysu.edu.cn(L. Chen), liuchangj@mail.sysu.edu.cn(C. Liu).}

\date{}
\maketitle
\begin{abstract}
This paper investigates the multiplicity and the number of limit cycles for planar piecewise linear system divided into two regions by a straight line and each linear subsystem has a node. Through constructing Poincar\'e half maps and a successor function, and analyzing the properties of the successor function, we can derive that this system has at most two limit cycles, counting the multiplicities of limit cycles.
\end{abstract}

\vspace{1mm}
\section{Introduction and statement of the main results}
Piecewise smooth systems are ubiquitous in diverse industries encompassing electrical circuits, automation, robotics control, chemical process management, biomedical engineering, and numerous other domains.

Hilbert's 16th problem's second component, central to the planar differential system's qualitative theory, explores the uniform upper bound $H(n)$ on the number of limit cycles for the polynomial system of degree $n.$ This problem naturally applies to planar piecewise smooth systems, though the latter poses greater challenges than smooth systems. The field of limit cycles of the piecewise smooth system was first explored in 1966 by Andronov et al. \cite{Andronov1966}.

Despite extensive research, the maximum number of limit cycles remains unknown, even for the most basic planar piecewise linear (PWL) system divided into two regions by a straight line. Two primary challenges hinder progress: one is the time at which the solution (we define the orbits arriving at the discontinuity line by using the Filippov convex approach, see \cite{Kuznetsov2003}) starts at the switching line and returning to the switching line again is hard to explain explicitly. The other is that, even after simplification, the system contains at least five parameters, which makes analysis quite challenging. According to \cite{Bernardo2008}, the limit cycles can be crossing or sliding. In our study, we consider the crossing limit cycles of PWL systems divided into two regions by a straight line. The following limit cycles, unless otherwise noted, are crossing limit cycles.

The general PWL system with a straight line of separation can be denoted as
\begin{equation}
\label{GE}
\mathbf{x}'
=\left\{
\begin{array}{l}
\mathbf{A}^-\mathbf{x}+\mathbf{b}^-,~~~~~x<0, \vspace{0.2cm}\\
\mathbf{A}^+\mathbf{x}+\mathbf{b}^+,~~~~~x>0,
\end{array}
\right.
\end{equation}
where $\mathbf{x}=(x~y)^T,$ $\mathbf{A}^\pm=(a_{ij}^\pm)$ are $2\times2$ constant matrix, and $\mathbf{b}^\pm=(b_1^\pm,b_2^\pm)^T$ are constant vectors of $\mathbb{R}^2$. The separation line is $x=0$, and for convenience, we call the left (resp. right) subsystem for $x<0$ (resp. $x>0$).

When system \eqref{GE} is continuous, i.e., $a_{12}^-=a_{12}^+,$ $a_{22}^-=a_{22}^+$ and $\mathbf{b}^-=\mathbf{b}^+.$ According to Lum and Chua's conjecture in 1991 \cite{Lum1991}, they assert that the continuous system \eqref{GE} can have no more than one limit cycle. Subsequently, Freire et al. provided theoretical validation of this hypothesis in 1998 \cite{Freire-Ponce-1998}.

When system \eqref{GE} relaxes the assumption of continuity, the known conclusions are about the refracting system (the system without sliding set, i.e., $b = 0$). In 2013 \cite{Freire2013}, Freire et al. proposed a conjecture that the refracting system has at most one limit cycle. This conjecture has been subsequently verified in \cite{Medrado2015}, with further corroboration in \cite{Li2021} and related works \cite{HUAN2013, HUAN2014, WANG2019, WANG20191, Ponce2018, Litao2020, Carmona2}.
However, for discontinuous system \eqref{GE} with the sliding set (i.e., $b\neq0$), most of the known results \cite{Llibre2011, Buzzi2013, HUAN2013, HUAN2014, WANG2019, WANG20191, Freire2014, Pessoa2022} are demonstrated with the lower bounds of the number of limit cycles. The conclusions are obtained by using a case-by-case (classified based on the singularity types of the subsystems) approach. With the exception of the case where system \eqref{GE} has a center, it is easy to see that system \eqref{GE} is classified into six cases, and known results are shown in Table \ref{table}. In the table, the vertical (resp. horizontal) axis denotes the type of singularity for the left (resp. right) subsystem, where $F$, $S$ and $N$ represent focus, saddle, and node.
\begin{table}[!htbp]
\centering
\begin{tabular}{p{3cm}p{3cm}p{3cm}p{2cm}}
\toprule
& F & S & N \\
 \midrule
F & 3 & 3 & 3 \\
S &   & 2 & 2 \\
N &   &   & 2 \\
\bottomrule
\end{tabular}
\caption{The lower bounds for the maximum number of limit cycles of system \eqref{GE}.}
\label{table}
\end{table}

It is natural to ask whether these lower bounds in  Table \ref{table} also serve as upper bounds. Several studies have provided upper bounds for specific cases, but a comprehensive understanding of different cases in Table \ref{table} is lacking. For example, when one of the subsystems contains a center or possesses two nonzero eigenvalues with opposing signs, it has been shown that there are at most two limit cycles in \cite{Llibre2018, Li2017}. Giannakopoulos et al. \cite{Giannakopoulos2001, Giannakopoulos2002} have proven that system with $\mathbb{Z}_2$-symmetry has at most two limit cycles. Furthermore, Llibre et al. \cite{Llibre2015, Llibre20151} have demonstrated that if the equilibrium of one subsystem lies on the discontinuity line, the upper bound remains two. If neither subsystem has an equilibrium, the upper bound reduces to one \cite{Llibre2017}.

Recently, Carmona et al. \cite{Carmona3} established that the number of limit cycles of system \eqref{GE} is upper bounded by $L^*,$ a natural number not exceeding $8$. Their approach uses an integral characterization of Poincar\'e half-maps associated with system \eqref{GE}, the foundational theories and applications detailed in \cite{Carmona0, Carmona1, Carmona2}. Obviously, $8$ would not be a sharp upper bound on the maximum number of limit cycles of system \eqref{GE}, so the upper bound on the maximum number of limit cycles of system \eqref{GE} remains an interesting open problem.

Our paper aims to investigate the upper bound on the maximum number of limit cycles of system \eqref{GE} with two nodes. There are three different  nodes:
a diagonalizable node with different eigenvalues ($N$); a star node, i.e., a diagonalizable node with equal eigenvalues ($N^*$); and an improper node, i.e., a non-diagonalizable node ($IN$). The cases $N^*N$, $N^*N^*$ or $N^*IN$ have no limit cycle obviously; the case of $ININ$ has been solved in \cite{Chen2024}, which proves that the system has at most one limit cycle; the cases $NIN$ and $NN$ have not been solved. Here, we only consider system \eqref{GE} with $NN$ case, and the $NIN$ case will be studied later. Our main result is as follows.

\begin{Theorem}\label{main th}
System \eqref{GE} with $NN$ type has at most two limit cycles, counting the multiplicities of limit cycles.
\end{Theorem}

An outline of this paper is as follows. Section \ref{s2} introduces preliminary concepts, including the canonical system form, definitions of Poincar\'e half maps and a successor function, as well as other essential knowledge. Subsequently, in Section \ref{s3} and Section \ref{s4}, we employ theoretical analysis to prove our main conclusions. In Section \ref{s5}, a discussion of the obtained results and potential strategies for their future promotion is provided.

\section{Preliminaries}
\label{s2}
According to reference \cite{Freire-Ponce-SIAM2012},  system \eqref{GE} and the following system are topologically equivalent when considering the number of crossing limit cycles, subject to the necessary condition that $a_{12}^-\cdot a_{12}^+>0$.
\begin{equation}
\label{NFE}
\left(
\begin{array}{cccc}
x'\\
y'\\
\end{array}
\right)
=
\left\{
\begin{array}{l}
\left(
\begin{array}{cccc}
2\gamma_{L}& -1\\
\gamma_{L}^2-m_{L}^2& 0\\
\end{array}
\right)\left(
\begin{array}{cccc}
x\\
y\\
\end{array}
\right)
-
\left(
\begin{array}{cccc}
0\\
\alpha_{L}\\
\end{array}
\right)~~~~x<0, \vspace{0.2cm}\\
\left(
\begin{array}{cccc}
2\gamma_{R}& -1\\
\gamma_{R}^2-m_{R}^2& 0\\
\end{array}
\right)\left(
\begin{array}{cccc}
x\\
y\\
\end{array}
\right)
-
\left(
\begin{array}{cccc}
-b\\
\alpha_{R}\\
\end{array}
\right)~~~~x>0,
\end{array}
\right.
\end{equation}
where $\gamma_{L,R},\alpha_{L,R},b\in\mathbb{R}$, and the modal parameters $m_{L,R}\in\{\mathrm{i}, 0, 1\}$ with $\mathrm{i}^2=-1.$
\begin{Remark}
System \eqref{NFE} has a focus for $m_{L,R}=\mathrm{i}$; has a node for $m_{L,R}=1$ and $|\gamma_{L,R}|>1$; has a saddle for $m_{L,R}=1$ and $|\gamma_{L,R}|<1$; and has an improper node for $m_{L,R}=0$ and $\gamma_L\gamma_R\ne 0$.
\end{Remark}

Thus, for system \eqref{NFE} with two nodes, we have $m_{L,R}=1$ and $|\gamma_{L,R}|>1$, system \eqref{NFE} will be written as
\begin{equation}
\label{NNE}
\left(
\begin{array}{cccc}
x'\\
y'\\
\end{array}
\right)
=
\left\{
\begin{array}{l}
\left(
\begin{array}{cccc}
2\gamma_{L}& -1\\
\gamma_{L}^2-1& 0\\
\end{array}
\right)\left(
\begin{array}{cccc}
x\\
y\\
\end{array}
\right)
-
\left(
\begin{array}{cccc}
0\\
\alpha_{L}\\
\end{array}
\right)~~~~x<0, \vspace{0.2cm}\\
\left(
\begin{array}{cccc}
2\gamma_{R}& -1\\
\gamma_{R}^2-1& 0\\
\end{array}
\right)\left(
\begin{array}{cccc}
x\\
y\\
\end{array}
\right)
-
\left(
\begin{array}{cccc}
-b\\
\alpha_{R}\\
\end{array}
\right)~~~~x>0.
\end{array}
\right.
\end{equation}

The equilibria of the left and right subsystems of system \eqref{NNE} are derived as
$$E_L=\left(\frac{\alpha_L}{\gamma_L^2-1}, \frac{2 \alpha_L \gamma_L}{\gamma_L^2-1}\right),~~E_R=\left(\frac{\alpha_R}{\gamma_R^2-1}, \frac{2\alpha_R\gamma_R}{\gamma_R^2-1}+b\right),$$ respectively. Upon straightforward calculation, the invariant straight lines for each subsystem are identified as
$$l_L^-:~y_L^-=(\gamma_L-1)x+\frac{\alpha_{L}}{\gamma_L-1},~l_L^+:~y_L^+=(\gamma_L+1)x+\frac{\alpha_{L}}{\gamma_L+1}$$
for the left subsystem, and
$$l_R^-:~y_R^-=(\gamma_R-1)x+\frac{\alpha_{R}}{\gamma_R-1}+b,~l_R^+:~y_R^+=(\gamma_R+1)x+\frac{\alpha_{R}}{\gamma_R+1}+b$$
for the right subsystem.

Subsequently, by neglecting the subscripts in \eqref{NNE}, we proceed to analyze the simplified system
\begin{equation}
\label{subs}
\left(
\begin{array}{cccc}
x'\\
y'\\
\end{array}
\right)
=\left(
\begin{array}{cccc}
2\gamma& -1\\
\gamma^2-1& 0\\
\end{array}
\right)\left(
\begin{array}{cccc}
x\\
y\\
\end{array}
\right)
-
\left(
\begin{array}{cccc}
-b\\
\alpha\\
\end{array}
\right).
\end{equation}

The solution of system \eqref{subs} starting from $\boldsymbol{x_0}=(x(0),y(0))^T$ is given by
\begin{equation}
\label{subss}
\left(
\begin{array}{cccc}
x(t)\\
y(t)\\
\end{array}
\right)
=
\left(
\begin{array}{cccc}
-\frac{\gamma-1}{2}&\frac{\gamma+1}{2}\\
-\frac{\gamma^2-1}{2}&\frac{\gamma^2-1}{2}\\
\end{array}
\right)
\left(
\begin{array}{cccc}
e^{(\gamma-1)t}&-\frac{e^{(\gamma-1)t}}{\gamma-1}\\
e^{(\gamma+1)t}&\frac{e^{(\gamma+1)t}}{\gamma+1}\\
\end{array}
\right)
\left(
\begin{array}{cccc}
x(0)-\bar{x}\\
y(0)-\bar{y}\\
\end{array}
\right)
+
\left(
\begin{array}{cccc}
\bar{x}\\
\bar{y}\\
\end{array}
\right),
\end{equation}
where $(\bar{x},\bar{y})^T$ is the equilibrium of system \eqref{subs}.

It follows from Proposition $2.2$ in \cite{HUAN2014} that a necessary condition for the existence of a periodic orbit of system \eqref{NNE} is stated. Here, we present both the proposition and a simple proof to facilitate subsequent usage.

\begin{Proposition}
\label{alar}
The necessary condition for the existence of a periodic orbit of system \eqref{NNE} is $\alpha_L>0>\alpha_R.$
\end{Proposition}
\noindent\bpf
Consider the vector fields  $$\boldsymbol{F}^\pm(\boldsymbol{x})=(F_1^\pm(x,y),F_2^\pm(x,y))^T$$ of the left and right subsystems, respectively.
Then
\begin{equation}
\label{F1}
F_1^-(0,y)=-y,~~F_1^+(0,y)=-y+b.
\end{equation}

Analysis of these fields along the $y$-axis reveals that, under their influence, orbits originating from $(0,y_0^-)$ (resp. $(0,y_0^+)$) with $y_0^->0$ (resp. $y_0^+<b$), will enter the region $S^-=\{(x,y):x<0\}$ (resp. $S^+=\{(x,y):x>0\}$) and then eventually intersect the $y$-axis again at $(0,y_1^-)$ (resp. $(0,y_1^+)$) with $y_1^-<0$ (resp. $y_1^+>b$), in finite time $t^->0$ (resp. $t^+>0$).

For a periodic orbit to exist, it must intersect the $y$-axis twice, at
$(0,y_0)$ and $(0,y_1),$ with conditions: $y_0>y_1,~y_0>\max\{0,b\},$ and $ y_1<\min\{0,b\}.$

Substituting the initial state $(x(0),y(0))=(0,y_0)$ into the equation \eqref{subss}, we let $x(t)=0,$  and get
\begin{equation*}
y_0=\frac{\alpha_L\big((\gamma_L+1)e^{(\gamma_L-1)t}-(\gamma_L-1)e^{(\gamma_L+1)t}-2\big)}{(\gamma_L^2-1)\big(e^{(\gamma_L-1)t}-e^{(\gamma_L+1)t}\big)},
\end{equation*}
where $(\gamma_{L}+1)e^{(\gamma_{L}-1)t}-(\gamma_{L}-1)e^{(\gamma_{L}+1)t}-2<0$ and $e^{(\gamma_{L}-1)t}-e^{(\gamma_{L}+1)t}<0$ for any $t>0.$

Thus, it can be shown that the necessary conditions for such intersections are $\alpha_L>0$ for the left subsystem and $\alpha_R<0$ for the right subsystem.\epf

\vspace{2mm}
Henceforth, we maintain the following assumption in this paper.

\begin{Assumption}
\label{parameter-condition}
   The parameters of system \eqref{NNE} satisfy: $|\gamma_{L,R}|>1$ and $\alpha_L>0>\alpha_R.$
\end{Assumption}

\subsection{The left Poincar\'e half map \texorpdfstring{$P_L(\cdot)$}{}}
Based on the analysis of Proposition \ref{alar}, it is established that the orbit originating from point $(0,y_0)$ with $y_0>0,$ under the flow of the left linear subsystem, will enter the left region $S^-$. If this orbit intersects the $y$-axis again at $(0,y_1)$ with $y_1<0$ after a positive time $t,$ then we can define a left half Poincar\'e map $P_L$ as
$$P_L(0)=0,~P_L(y_0)=y_1,~y_0>0.$$

The parametric representation of the Poincar\'e half map $P_L$ is given by
\begin{equation}
\label{ly0}
y_0^L=\frac{\alpha_L\big((\gamma_L+1)e^{(\gamma_L-1)t}-(\gamma_L-1)e^{(\gamma_L+1)t}-2\big)}{(\gamma_L^2-1)\big(e^{(\gamma_L-1)t}-e^{(\gamma_L+1)t}\big)},
\end{equation}
\begin{equation}
\label{ly1}
y_1^L=\frac{\alpha_L\big((\gamma_L-1)e^{(\gamma_L-1)t}-(\gamma_L+1)e^{(\gamma_L+1)t}+2e^{2\gamma_Lt}\big)}{(\gamma_L^2-1)\big(e^{(\gamma_L-1)t}-e^{(\gamma_L+1)t}\big)},~~t>0.
\end{equation}

\vspace{2mm}
By Proposition $2.3$ in \cite{HUAN2014}, some properties of $P_L(\cdot)$ are stated below.
\begin{Proposition}
\label{pl}
Analyzing the Poincar\'e half map $P_L$, we have the following properties.
\begin{description}
\item $(i)$ $y_0$ is increasing and $y_1$ is decreasing with respect to $t.$
\item $(ii)$ $P_L$ has $y_0=\frac{\alpha_L}{\gamma_L+1}$ as an asymptote for $\gamma_L>0$ and has $y_1=\frac{\alpha_L}{\gamma_L-1}$ as an asymptote for $\gamma_L<0$.
\item $(iii)$ When $\gamma_L>0,$ the domain of definition for $P_L$ is $(0,\frac{\alpha_L}{\gamma_L+1}),$ $P_L$ is decreasing and concave down with respect to $y_0.$ When $\gamma_L<0,$ the domain of definition for $P_L$ is $(0,+\infty),$ $P_L$ is decreasing and concave up with respect to $y_0.$
\item $(iv)$ We define $P_L(0)=0,$ then $P_L$ is continuous at $y_0=0.$ What's more, the first four derivatives of $P_L$ at $y_0=0$ are
$$P_L'(0)=-1,~~P_L''(0)=-\frac{8\gamma_L}{3\alpha_L},~P_L'''(0)=-\frac{32\gamma_L^2}{3\alpha_L^2},~
P_L^{IV}(0)=-\frac{32\gamma_L(23\gamma_L^2+9)}{9\alpha_L^3}.$$
\end{description}
\end{Proposition}

\subsection{The right Poincar\'e half map \texorpdfstring{$P_R^{-1}(\cdot)$}{}}
Similarly, the orbit from point $(0,y_1)$ with $y_1<b,$ under the flow of the right linear subsystem, will enter the right region $S^+$. If this orbit intersects the $y$-axis again at $(0,y_0)$ with $y_0>b$ after a positive time $t,$ then we can define a right Poincar\'e half map $P_R$ as
$$P_R(b;b)=b,~P_R(y_1;b)=y_0,~y_1<b.$$

The inverse mapping  $P_R^{-1}(y_0;b)$ of $P_R(y_1;b)$ is naturally defined as
$$P_R^{-1}(b;b)=b,~P_R^{-1}(y_0)=y_1,~y_0>b.$$

The parametric representation of the Poincar\'e half map $P_R^{-1}$ is given by
\begin{equation}
\label{ry0}
y_0^R=\frac{\alpha_R\big((\gamma_R+1)e^{(\gamma_R-1)s}-(\gamma_R-1)e^{(\gamma_R+1)s}-2\big)}{(\gamma_R^2-1)\big(e^{(\gamma_R-1)s}-e^{(\gamma_R+1)s}\big)}+b,
\end{equation}
\begin{equation}
\label{ry1}
y_1^R=\frac{\alpha_R\big((\gamma_R-1)e^{(\gamma_R-1)s}-(\gamma_R+1)e^{(\gamma_R+1)s}+2e^{2\gamma_Rs}\big)}{(\gamma_R^2-1)\big(e^{(\gamma_R-1)s}-e^{(\gamma_R+1)s}\big)}+b,~~s<0.
\end{equation}

\vspace{2mm}
By Proposition $2.4$ in \cite{HUAN2014}, some properties of $P_R^{-1}(\cdot)$ are stated below.
\begin{Proposition}
\label{pr}
Analyzing the Poincar\'e half map $P_R^{-1}$, we have the following properties.
\begin{description}
\item $(i)$ $y_0(s)$ is increasing and $y_1(s)$ is decreasing with respect to $s.$
\item $(ii)$ $P_R^{-1}$ has $y_1=\frac{\alpha_R}{\gamma_R+1}+b$ as an asymptote for $\gamma_R>0$ and has $y_0=\frac{\alpha_R}{\gamma_R-1}+b$ as an asymptote for $\gamma_R<0$.
\item $(iii)$ When $\gamma_R>0,$ the domain of definition for $P_R^{-1}$ is
$(b,+\infty),$ $P_R^{-1}$ is decreasing and concave up with respect to $y_0.$ When $\gamma_R<0,$ the domain of definition for $P_R^{-1}$ is
$(b,\frac{\alpha_R}{\gamma_R-1}+b),$ $P_R^{-1}$ is decreasing and concave down with respect to $y_0.$
\item $(iv)$ We define $P_R^{-1}(b;b)=b,$ then $P_R^{-1}$ is continuous at $y_0=b.$ What's more, the first four derivatives of $P_R^{-1}$ at $y_0=b$ are
$$\big(P_R^{-1}\big)'(b;b)=-1,~~\big(P_R^{-1}\big)''(b;b)=-\frac{8\gamma_R}{3\alpha_R},$$
$$\big(P_R^{-1}\big)'''(b;b)=-\frac{32\gamma_R^2}{3\alpha_R^2},~
\big(P_R^{-1}\big)^{IV}(b;b)=-\frac{32\gamma_R(23\gamma_R^2+9)}{9\alpha_R^3}.$$
\end{description}
\end{Proposition}

\subsection{The successor function  \texorpdfstring{$d(\cdot;b)$}{}  }\label{successor}
The successor function $d(\cdot;b)$ for a specified $b$ is defined by utilizing the left Poincar\'e half map $P_L(\cdot)$ and the inverse of the right Poincar\'e half map $P_R^{-1}(\cdot)$ as follows.
\begin{equation}\label{dy0}
d(y_0;b)=P_R^{-1}(y_0;b)-P_L(y_0),~~y_0\in\big(y_0^m,y_0^M\big),
\end{equation}
the values $y_0^m$ and $y_0^M$ are derived from  Propositions \ref{pl} and \ref{pr}.
\begin{Proposition}
\label{y0M}
    $y_0^m=\max\{0,b\}$ and the value of $y_0^M$ has four cases as follows.
    \begin{description}
    \item $(i)$  $y_0^M=\frac{\alpha_L}{\gamma_L+1}$ for $\gamma_L,\gamma_R>0.$
    \item $(ii)$ $y_0^M=\min\big\{\frac{\alpha_L}{\gamma_L+1},\frac{\alpha_R}{\gamma_R-1}+b\big\}$ for $\gamma_L>0>\gamma_R,$ where $\frac{\alpha_R}{1-\gamma_R}<b<\frac{\alpha_L}{\gamma_L+1}.$
    \item $(iii)$ $y_0^M=+\infty$ for $\gamma_L<0<\gamma_R.$
    \item $(iv)$ $y_0^M=\frac{\alpha_R}{\gamma_R-1}+b$ for $\gamma_L,\gamma_R<0,$      where   $b>\frac{\alpha_R}{1-\gamma_R}.$    \end{description}
\end{Proposition}

When $b=0$, the Taylor series expansion of $P_L(y_0)$ and $P_R^{-1}(y_0;0)$ at $y_0=0$ are
$$P_L(y_0)=-y_0-\frac{4\gamma_L}{3\alpha_L}y_0^2-\frac{16\gamma_L^2}{9\alpha_L^2}y_0^3-
\frac{4\gamma_L(23\gamma_L^2+9)}{27\alpha_L^3}y_0^4+o(y_0^4),$$
and
$$P_R^{-1}(y_0;0)=-y_0-\frac{4\gamma_R}{3\alpha_R}y_0^2-\frac{16\gamma_R^2}{9\alpha_R^2}y_0^3-
\frac{4\gamma_R(23\gamma_R^2+9)}{27\alpha_R^3}y_0^4+o(y_0^4),$$
respectively. Leading to the Taylor series expansion of $d(y_0;0)$ at $y_0=0$ as
\begin{equation}\label{dy00}
d(y_0;0)=\Bigg(\frac{4\gamma_L}{3\alpha_L}-\frac{4\gamma_R}{3\alpha_R}\Bigg)y_0^2+
\frac{16}{9}\Bigg(\frac{\gamma_L^2}{\alpha_L^2}-\frac{\gamma_R^2}{\alpha_R^2}\Bigg)y_0^3+\frac{4}{27}
\Bigg(\frac{\gamma_L(23\gamma_L^2+9)}{\alpha_L^3}-\frac{\gamma_R(23\gamma_R^2+9)}{\alpha_R^3}\Bigg)y_0^4
+o(y_0^4).
\end{equation}

For fixed $b$, the (multiple) root $y_0^*$ of $d(y_0;b)$ corresponds to a (multiple) limit cycle of system \eqref{NNE}, and furthermore, the sign of $d'(y_0^*;b)$ determines the stability of the limit cycle, here '$'$' denotes the first-order partial derivative of $d(y_0;b)$ with respect to $y_0.$ For convenience, we will sometimes say that the root $y_0^*$ represents a limit cycle.

In the following, we will give some properties on $d(y_0;b).$ These properties have been described in literature \cite{Freire-Ponce-SIAM2012}, and here, for the reader's convenience, we give a simple proof.
\begin{Proposition}
\label{propdy0}
According to the definition of $d(y_0;b),$ the function $d(\cdot;b)$ has the following conclusions.
\begin{description}
 \item $(i)$ $d(y_0;b)=d(y_0-b;0)+b$ and $\frac{\partial d(y_0;b)}{\partial b}>0,$ where $y_0\in\big(y_0^m,y_0^M\big).$
 \item $(ii)$ $sign(d(y_0;b))=sign(P_R^{-1}(y_0;b)-P_L(y_0))=sign(P(y_0;b)-y_0),$ where $P(y_0;b)=P_R\circ P_L(y_0).$
 \item $(iii)$ Suppose $y_0^*$ is a limit cycle of system \eqref{NNE} for a given $b$, where $y_0^*\in\big(y_0^m,y_0^M\big).$ $y_0^*$ is stable for $d'(y_0^*;b)<0$ and unstable for $d'(y_0^*;b)>0.$ If $d'(y_0^*;b)=0,$ then $y_0^*$ is internally stable and externally unstable for $d''(y_0^*;b)>0$ and $y_0^*$ is internally unstable and externally stable for $d''(y_0^*;b)<0.$
\end{description}
\end{Proposition}
\noindent\bpf
$(i)$ Define $Graph(\cdot)$ is the set of all points of the graph of a function. Taking $(y_0,y_1)\in Graph\big(P_R^{-1}(\cdot;b)\big),$ then $y_1=P_R^{-1}(y_0;b).$
By the transformation $(x,y-b,t)\rightarrow (x,y,t),$ the parameter $b$ in the right subsystem can be eliminated, thus
$$(y_0-b,y_1-b)\in Graph(P_R(\cdot;0)),$$
we have $y_1-b=P_R^{-1}(y_0-b;0).$ i.e.,
$$P_R^{-1}(y_0;b)=P_R^{-1}(y_0-b;0)+b.$$
Thus,
$$d(y_0;b)=P_R^{-1}(y_0;b)-P_L(y_0)=P_R^{-1}(y_0-b;0)
+b-P_L(y_0)=d(y_0-b;0)+b,$$
and
$$\frac{\partial d(y_0;b)}{\partial b}=\frac{\partial d(y_0-b;0)}{\partial b}+1=-\frac{\partial P_R^{-1}(y_0-b;0)}{\partial (y_0-b)}+1>0.$$

$(ii)$ Since $d(y_0;b)=P_R^{-1}(y_0;b)-P_L(y_0)=P_R^{-1}(y_0;b)-P_R^{-1}(P(y_0);b)$ and $P_R^{-1}(y_0;b)$ is monotonically decreasing with respect to $y_0$, we have $sign(P_R^{-1}(y_0;b)-P_L(y_0))=sign(P(y_0;b)-y_0)$

\begin{figure}[htpb!]
\centering  
\includegraphics[width=0.3\textwidth]{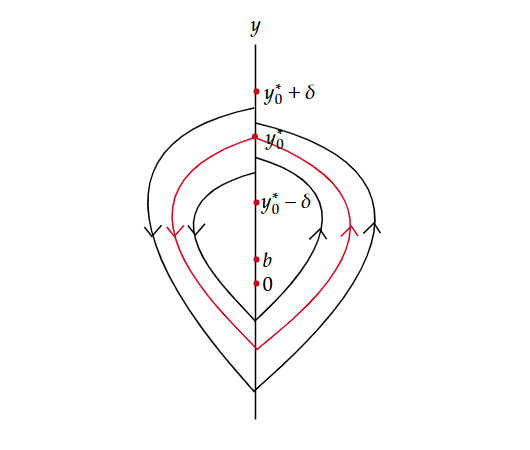}
\caption{A stable limit cycle $L_{y_0^*}$ throughing $(0,y_0^*)^T$ of system \eqref{NNE} for fixed $b$}
\label{fig1}
\end{figure}

$(iii)$ We first assume that $d(y_0^*;b)=0$ and $d'(y_0^*;b)<0,$ then there exists a $\delta,$ where $0<\delta\ll1,$ such that $d(y_0;b)$ is monotonically decreasing on $\big(y_0^*-\delta, y_0^*+\delta\big)$. Since $d(y_0^*;b)=0,$ we have $d(y_0;b)>0$ for $y_0\in\big(y_0^*-\delta, y_0^*\big)$ and $d(y_0;b)<0$ for $y_0\in\big(y_0^*, y_0^*+\delta\big)$. Then, by the definition of $d(y_0;b)$, we have $P_R^{-1}(y_0;b)>P_L(y_0)$ when $y_0\in\big(y_0^*-\delta, y_0^*\big)$ and $P_R^{-1}(y_0;b)<P_L(y_0)$ when $y_0\in\big(y_0^*, y_0^*+\delta\big),$ which means that $y_0^*$ is stable. The schematic diagram is shown as Fig. \ref{fig1}.
Similarly, we can obtain that $y_0^*$ is unstable if $d'(y_0^*;b^*)>0.$

When $d(y_0^*;b)=d'(y_0^*;b)=0$ and $d''(y_0^*;b)>0,$ there exist a $\delta,$ where $0<\delta\ll1,$  we have $d'(y_0;b)<0$ for $y_0\in\big(y_0^*-\delta, y_0^*\big)$ and $d'(y_0;b)>0$ for $y_0\in\big(y_0^*, y_0^*+\delta\big)$. That is, $y_0^*$ is an internally stable and externally unstable semi-stable limit cycle. Similarly, we can obtain that $y_0^*$ is internally unstable and externally stable if $d'(y_0^*;b)=0$ and $d''(y_0^*;b)<0.$
\epf

\section{The multiplicity and number of limit cycles of system  \texorpdfstring{\eqref{NNE}}{}}
\label{s3}
Firstly, we present three crucial parameters
\begin{equation}\label{delta123}
  \Delta_1=\frac{\gamma_L}{\alpha_L}-\frac{\gamma_R}{\alpha_R},~
  \Delta_2=\frac{\alpha_R}{\gamma_R-1}-\frac{\alpha_L}{\gamma_L+1},~
\Delta_3=\frac{\alpha_R}{\gamma_R+1}-\frac{\alpha_L}{\gamma_L-1}.
\end{equation}

When $b=0,$ system \eqref{NNE} is a refracting system with two nodes. The number and the existence of limit cycles of such a system have been obtained in Theorem $3.1$ of \cite{HUAN2014}.

\begin{Theorem}\label{beq0lc}
Suppose $b=0$ and Assumption \ref{parameter-condition} holds, then we have the following statements.
\begin{description}
\item $(i)$ If $\gamma_L\gamma_R>0,$ system \eqref{NNE} has no limit cycle.
\item $(ii)$ If $\Delta_1=\Delta_2=0,$ the origin is a local center.
\item $(iii)$ If $\gamma_L>0>\gamma_R,$ we have the following statements.
    \item $(a)$ When $\Delta_1\cdot\Delta_2\geq0$ and $\Delta_1+\Delta_2\neq0,$ system \eqref{NNE} has no limit cycle.
    \item $(b)$ When $\Delta_1\cdot\Delta_2<0,$ system \eqref{NNE} has a unique stable (resp. unstable) limit cycle if $\Delta_2<0$ (resp. $\Delta_2>0$).
        \item $(iv)$ If $\gamma_L<0<\gamma_R,$ we have the following statements.
    \item $(a)$ When $\Delta_1\cdot\Delta_3\geq0$ and $\Delta_1+\Delta_3\neq0,$ system \eqref{NNE} has no limit cycle.
    \item $(b)$ When $\Delta_1\cdot\Delta_3<0,$ system \eqref{NNE} has a unique stable (resp. unstable) limit cycle if $\Delta_3<0$ (resp. $\Delta_3>0$).
\end{description}
\end{Theorem}

In addition, when $\gamma_L\gamma_R>0$ and $b\neq0,$ the literature \cite{HUAN2014} (Theorem $3.2$) establishes that system \eqref{NNE}
has at most one limit cycle for system \eqref{NNE}.
\begin{Theorem}\label{rlrrg0}
Suppose $\gamma_L\gamma_R>0$ and Assumption \ref{parameter-condition} holds, then we have the following statements.
\begin{description}
\item $(i)$ If $b\gamma_R>0,$ no periodic orbit exists.
\item $(ii)$ If $b\gamma_R<0,$ a unique periodic orbit exists that is stable when $b>0$ and unstable when $b<0.$
\end{description}
\end{Theorem}

Subsequently, we focus on the remaining case $\gamma_L\gamma_R<0$ and $b\neq0.$ Assuming $\gamma_L>0>\gamma_R$ without loss of generality. For $\gamma_L<0<\gamma_R,$ the conclusions mirror those of $\gamma_L>0>\gamma_R$ through a symmetry transformation
\begin{equation*}
(x,y,t,\gamma_L,\gamma_R,\alpha_L,\alpha_R,b)
\to(x,-y,-t,-\gamma_L,-\gamma_R,\alpha_L,\alpha_R,-b).
\end{equation*}

\begin{Theorem}
\label{pseudohopf}
Suppose $\gamma_L>0>\gamma_R$ and Assumption \ref{parameter-condition} holds. Then, for each $|b|<\varepsilon$ that is sufficiently small and satisfies $b\Delta_1<0,$ system \eqref{NNE} can generate a unique limit cycle originating from the origin. If $\Delta_1<0,$ the limit cycle is stable, if $\Delta_1>0,$ it is unstable. Additionally, when $\Delta_1=0,$ we need to determine the sign of $\Delta_2,$ the limit cycle is stable for $\Delta_2<0$ and unstable for $\Delta_2>0.$
\end{Theorem}
\noindent\bpf
For $\varepsilon>0,$ we define a Poincar\'e map $P:(0,\varepsilon)\to(0,\varepsilon)$ with $P(y_0,b)=P_R\circ P_L(y_0).$ Utilizing the properties of $P_L$ and $P_R$ and Propositions \ref{pl} and \ref{pr}, the Taylor expansion of $P(y_0,b)$ at $(0,0)$ is derived as
$$P(y_0,b)=y_0+\frac{4}{3}\Big(\frac{\gamma_{L}}{\alpha_L}-\frac{\gamma_{R}}{\alpha_R}\Big)y_0^2+o\big(y_0^2\big)+g_0(b)+\sum_{k=1}g_k(b)y_0^k,$$
where $g_0(b)=2b+O\big(|b|^2\big)$ and $g_k(b)=O\big(|b|\big)$ for $k\geq1.$

Defining $$P_1(y_0,b)=P(y_0,b)-y_0=\frac{4}{3}\Big(\frac{\gamma_{L}}{\alpha_L}-\frac{\gamma_{R}}{\alpha_R}\Big)y_0^2+o\big(y_0^2\big)+g_0(b)+\sum_{k=1}g_k(b)y_0^k,$$
we note that $P_1(0,0)=0$ and $\frac{\partial P_1}{\partial b}(0,0)=2.$
The implicit function theorem ensures the existence of a unique function
$$b=b(y_0)=\frac{2}{3}\Big(\frac{\gamma_{R}}{\alpha_R}-\frac{\gamma_{L}}{\alpha_L}\Big)y_0^2+o\big(y_0^2\big),$$
such that $P_1\big(y_0,b(y_0)\big)=0$ for any $y_0\in(0,\varepsilon),$ where $0<\varepsilon\ll1.$ Specifically, for sufficiently small $b$ with $b\Big(\frac{\gamma_{L}}{\alpha_L}-\frac{\gamma_{R}}{\alpha_R}\Big)=b\Delta_1<0,$ there exists a $y_0\in(0,\varepsilon)$ satisfying $P\big(y_0,b(y_0)\big)=y_0.$

Finally, consider the stability of the limit cycle. The partial derivative of $P$ with respect to $y_0$ is
$$\frac{\partial}{\partial y_0}P(y_0,b)=1+\frac{8}{3}\Big(\frac{\gamma_{L}}{\alpha_L}-\frac{\gamma_{R}}{\alpha_R}\Big)y_0+\sum_{k=1}g_k(b)y_0^{k-1}+o\big(y_0\big)$$
and for each $b=b(y_0),~y_0\in(0,\varepsilon),$ we have
$$\frac{\partial}{\partial y_0}P\big(y_0,b(y_0)\big)=1+\frac{8}{3}\Big(\frac{\gamma_{L}}{\alpha_L}-\frac{\gamma_{R}}{\alpha_R}\Big)y_0+o\big(y_0\big).$$

This is less than 1 for $\Delta_1<0$ and greater than 1 for $\Delta_1>0,$ so that when $\Delta_1<0,$ the limit cycle is stable, and when $\Delta_1>0,$ the limit cycle is unstable.

When $b\Delta_1>0,$ $P_1\big(y_0,b(y_0)\big)$ does not vanish, and therefore, there is no periodic orbit of system \eqref{NNE}.

When $\Delta_1=0,$ the Taylor expansion of $P(y_0,b)$ at $(0,0)$ is derived as
$$P(y_0,b)=y_0+\frac{4}{27}\Bigg(\frac{\gamma_{L}\big(23\gamma_L^2+9\big)}{\alpha_L^3}-\frac{\gamma_{R}\big(23\gamma_R^2+9\big)}{\alpha_R^3}\Bigg)y_0^4+o\big(y_0^4\big)+g_0(b)+\sum_{k=1}g_k(b)y_0^k.$$

Since
\begin{align*}
\frac{\gamma_{L}\big(23\gamma_L^2+9\big)}{\alpha_L^3}-\frac{\gamma_{R}\big(23\gamma_R^2+9\big)}{\alpha_R^3} &=\frac{23\gamma_R^3}{\alpha_R^3}+\frac{9\gamma_R}{\alpha_R\alpha_L^2}-\frac{23\gamma_R^3}{\alpha_R^3}-\frac{9\gamma_R}{\alpha_R^3}\\
&=\frac{9\gamma_R}{\alpha_R}\Big(\frac{1}{\alpha_L^2}-\frac{1}{\alpha_R^2}\Big)\\
&=\frac{9\gamma_R(\alpha_R-\alpha_L)(\alpha_L+\alpha_R)}{\alpha_R^3\alpha_L^2},
\end{align*}
and combine with $\Delta_2=\frac{\alpha_L+\alpha_R}{(\gamma_R-1)(\gamma_L+1)}$, the remaining proof is similarity to the case $\Delta_1\neq0$, we omit it for simplify.
\epf

\begin{Remark}
If $\Delta_1=\Delta_2=0,$ then according to statement $(II)$ of Theorem \ref{beq0lc}, we derive $y_0^M=\frac{\alpha_L}{\gamma_L+1}$ and $d(y_0;0)\equiv0$ for any $y_0\in\big(0,y_0^M\big)$. Notice that $\frac{\partial d(y_0;b)}{\partial b}>0,$ thus $d(y_0;b)\neq0$ for all $y_0\in\big(y_0^m,y_0^M\big)$ and $b\neq0$, that is, the system \eqref{NNE} has no periodic orbit.
\end{Remark}

\begin{Remark}
The aforementioned limit cycle in Theorem \ref{pseudohopf} formation can be succinctly articulated as transitioning from a state with $b\neq0$ to $b=0,$ resulting in a collision of two invisible tangency points. This collision represents a bifurcation, specifically the $II_2$ or pseudo-Hopf bifurcation, as described in \cite{Kuznetsov-Rinaldi-Gragnani-IJBC2003}.
\end{Remark}

In the following, we present a theorem that describes the variation of the roots of the function $d(y_0;b)$ with respect to the parameter $b.$ The core idea of this theorem is derived from Lemma $5$ in reference \cite{Carmona3}, which we have adapted and subsequently provided a concise proof to suit our specific proof context better.
\begin{Theorem}\label{dwithb}
According to the expression of $d(y_0;b),$ for a given $b=b^*$ and $y_0\in\big(y_0^m,y_0^M\big),$ we have the following conclusions.
\begin{description}
\item $(i)$ If $d(y_0;b^*)\equiv0,$ then $d(y_0;b)$ has no root for any $\mid b-b^*\mid$ sufficiently small.
\item$(ii)$ If $d(y_0;b^*)$ has a simple root $y_0=y_0^*,$ then for $\mid b-b^*\mid$ sufficiently small, $d(y_0;b)$ has a simple root $y_0=y_0^{1*}$, where $|y_0^{1*}-y_0^*|\ll 1$.  What's more, we have $y_0^{1*}>y_0^*$ when  $d'(y_0^*;b^*)<0$ and $y_0^{1*}<y_0^*$ when $d'(y_0^*;b^*)>0.$
\item$(iii)$ If $d(y_0;b^*)$ has a doubled root $y_0=y_0^*,$ then when $d''(y_0^*;b^*)<0$, for $\mid b-b^*\mid$ sufficiently small, in some neighborhood of $y_0^*$, $d(y_0;b)$ has no root for $b<b^*$ and has two roots $y_0=y_0^{1*}$ and $y_0=y_0^{2*}$ for $b>b^*,$ where $y_0^{1*}<y_0^{2*},~d'(y_0^{1*};b^*)>0$ and $d'(y_0^{2*};b^*)<0$.  When $d''(y_0^*;b^*)>0,$  in some neighborhood of $y_0^*$, $d(y_0;b)$ has no root for $b>b^*$ and has two roots $y_0=y_0^{1*}$ and $y_0=y_0^{2*}$ for $b<b^*,$ where $y_0^{1*}<y_0^{2*},~d'(y_0^{1*};b^*)<0$ and $d'(y_0^{2*};b^*)>0$.
\end{description}
\end{Theorem}
\noindent\bpf
$(i)$ If $d(y_0;b^*)\equiv0$ for any $y_0\in\big(y_0^m,y_0^M\big),$ then based on the fact that $\frac{\partial{d(y_0;b)}}{\partial{b}}>0$, we can deduce that $d(y_0;b^*)\neq0$ for any $y_0\in\big(y_0^m,y_0^M\big)$ and $b\neq b^*.$

$(ii)$ If $d(y_0;b^*)$ has a simple root $y_0=y_0^*,$ the Taylor series expansion of the function $d(y_0;b)$ at $(y_0^*,b^*)$ yields
$$d(y_0;b)=a(y_0-y_0^*)(1+ o(1))+m(b-b^*)(1+o(1)),$$
where $a=d'(y_0^*;b^*)$ and $m=\frac{\partial{d}}{\partial{b}}(y_0^*,b^*)>0.$

By the implicit function theorem, there exists a unique $y_0=y_0^{1*}(b)$ satisfying $d(y_0^{1*}(b);b)=0$ for $\mid b-b^*\mid$ sufficiently small, and the sign of $(y_0^{1*})'(b)$ is opposite to $a.$

$(iii)$ If $d(y_0;b^*)$ has a doubled root $y_0=y_0^*,$ the Taylor series expansion of the function $d(y_0;b)$ at $b=b^*$ and $y_0=y_0^*$ yields
$$d(y_0;b)=a(y_0-y_0^*)^2(1+o(1))+m(b-b^*)(1+o(1)).$$

Similar to the case in $(ii),$ if $a>0,$ for some $0<\varepsilon\ll1,$ two roots $y_0=y_0^{1*,2*}(b)$ exist for $b\in(b^*-\varepsilon,b^*),$ and no root exists for $b\in(b^*,b+\varepsilon).$ Conversely, if $a<0$, two roots $y_0=y_0^{1*,2*}(b)$ exist for $b\in(b^*,b+\varepsilon)$ and no root exists for $b\in(b^*-\varepsilon,b^*).$ The sign of $\big(y_0^{1*,2*}\big)'(b)$ can be determined directly.
\epf

\begin{Remark}\label{limitcycle}
Notice that the number of limit cycles corresponds to the zeros of $d(\cdot;b)$ on $\big(y_0^m,y_0^M\big),$ according to Theorem \ref{dwithb}, if system \eqref{NNE} has a closed orbit $L^*$ at $b=b^*$, then the following assertions hold for $|b-b^*|$ sufficiently small.
\begin{description}
\item $(i)$ If $L^*$ is a non-isolated closed orbit, then for $b\neq b^*$, system \eqref{NNE} has no closed orbit near $L^*$.
\item$(ii)$ If $L^*$ is a hyperbolic limit cycle, then system \eqref{NNE} has a unique limit cycle $L_b$ near $L^*$, and $L_b$ expands (resp. contracts) as $b$ increases if $L^*$ is a stable (resp. unstable) limit cycle.
\item$(iii)$ If $L^*$ is a limit cycle with a multiplicity of two, the internally unstable and externally stable limit cycle $L^*$ will split into two limit cycles near $L^*$ when $b$ increases and vanishes when $b$ decreases. The internally stable and externally unstable limit cycle $L^*$ will split into two limit cycles near $L^*$ when $b$ decreases and vanishes when $b$ increases.
\end{description}
\end{Remark}

\begin{Theorem}
\label{dy02}
Suppose that $\gamma_L>0>\gamma_R$ and Assumption \ref{parameter-condition} holds. Then, the multiplicity of limit cycles of system \eqref{NNE} is at most two. Furthermore, when system \eqref{NNE} has a limit cycle of multiplicity two,
i.e., there exists a $y_0^*\in\big(y_0^m,y_0^M\big),$ such that $d(y_0^*;b)=d'(y_0^*;b)=0,$ we can derive that $d''(y_0^*;b)<0$ for $\gamma_L>-\gamma_R$ and $d''(y_0^*;b)>0$ for $\gamma_L<-\gamma_R.$ That is, the semi-stable limit cycle is internally stable (resp. unstable) and externally unstable (resp. stable) when $\gamma_L+\gamma_R<0$ (resp. $\gamma_L+\gamma_R>0$).
\end{Theorem}

The proof of Theorem \ref{dy02} is given in Section \ref{s4} since it is complicated.

\vspace{2mm}

Next, we will analyze the number of limit cycles of system \eqref{NNE}. Here, we introduce three crucial values
$$b_m=\frac{\alpha_R}{1-\gamma_R},~~\bar{b}=\frac{\alpha_L}{\gamma_L+1}-\frac{\alpha_R}{1-\gamma_R},~~ b_M=\frac{\alpha_L}{\gamma_L+1}.$$

Firstly, we consider the case $b>0$. We assume that $\gamma_L>0>\gamma_R$ and Assumption \ref{parameter-condition} holds. Recall that the function $d(y_0;b)$ is well defined on the interval $(y_0^m, y_0^M)$, where $y_0^m=\max\{0, b\}=b$ and $y_0^M=\min\{\frac{\alpha_R}{\gamma_R-1}+b,\frac{\alpha_L}{\gamma_L+1}\}.$
Since $y_0^m$ must be less than $y_0^M$ for $d(y_0;b)$ to be well defined, it follows that $b<b_M$.

Now we consider $d(y_0;b)$ as a binary function in $(y_0, b)$, which is well defined on the region $\Omega=\{(y_0, b)\,|\, 0<b<b_M, b<y_0<y_0^M\}$. The graph of $\Omega$ has two possibilities.

\begin{figure}[htpb!]
\centering
\subfigure[$\frac{\alpha_L}{\gamma_L+1}\leq\frac{\alpha_R}{\gamma_R-1}$]{
\includegraphics[width=0.45\textwidth]{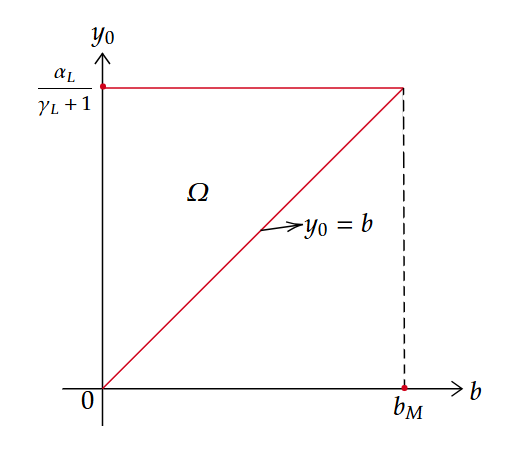}}
\subfigure[$\frac{\alpha_L}{\gamma_L+1}>\frac{\alpha_R}{\gamma_R-1}$]{
\includegraphics[width=0.45\textwidth]{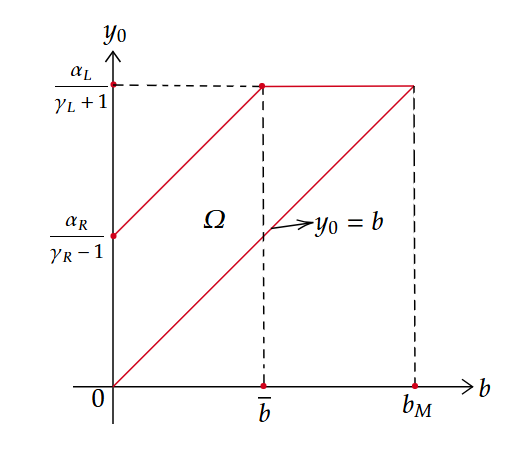}}
\caption{The graph of $\Omega$ for $b>0$}
\label{figy0b>0}
\end{figure}

$(1)$ When $\frac{\alpha_L}{\gamma_L+1}\leq\frac{\alpha_R}{\gamma_R-1},$ we have $y_0^M=\frac{\alpha_L}{\gamma_L+1}$ (see Fig. \ref{figy0b>0}$(a)$).

$(2)$ When $\frac{\alpha_L}{\gamma_L+1}>\frac{\alpha_R}{\gamma_R-1},$ there exists a unique $b=\bar{b},$ we have $y_0^M=\frac{\alpha_L}{\gamma_L+1}$ for $b\geq\bar{b}$ and $y_0^M=\frac{\alpha_R}{\gamma_R-1}+b$ for $b<\bar{b}$ (see Fig. \ref{figy0b>0}$(b)$).

Then we will study the sign of $d(y_0; b)$ on some boundaries of $\Omega$. According to Proposition \ref{propdy0} $(i),$ on the boundary $\{y_0=b\}$, we have
$$d(b;b)=d(b-b;0)+b=d(0;0)+b=b>0.$$

On the boundary $\{b=0\}$, it is known that $d(y_0;0)$ has at most one zero.

On the  boundary $\{y_0=\frac{\alpha_R}{\gamma_R-1}+b\}$, we have $\lim\limits_{y_0\to\frac{\alpha_R}{\gamma_R-1}+b}P_R^{-1}(y_0;b)=-\infty.$ Since $\lim\limits_{y_0\to\frac{\alpha_L}{\gamma_L+1}}P_L(y_0)=-\infty$ and $\frac{\alpha_R}{\gamma_R-1}+b<\frac{\alpha_L}{\gamma_L+1},$ it follows that $P_L\Big(\frac{\alpha_R}{\gamma_R-1}+b\Big)\in \mathbb R$. Therefore,
$$\lim\limits_{y_0\to \frac{\alpha_R}{\gamma_R-1}+b}d(y_0;b)=\lim\limits_{y_0\to \frac{\alpha_R}{\gamma_R-1}+b}P_R^{-1}(y_0;b)-P_L\Big(\frac{\alpha_R}{\gamma_R-1}+b\Big)=-\infty<0.$$

\begin{Theorem}\label{instable} Under the assumptions that $\gamma_L>0>\gamma_R$, $b>0$ and Assumption \ref{parameter-condition} holds. If $\gamma_L+\gamma_R>0$, then system \eqref{NNE} has at most two limit cycles.
\end{Theorem}
\noindent\bpf
If the conclusion does not hold,  then there exists a  $b_*\in(0,b_M)$ such that system \eqref{NNE} has at least three limit cycles. Without loss of generality, we assume that they are all hyperbolic for simplicity. If non-hyperbolic limit cycles exist, then their multiplicities must be two and they are inner unstable. By adjusting $b$ slightly (e.g., $b=b_*+\varepsilon$ with $0<\varepsilon\ll1$), any non-hyperbolic limit cycle would split into two hyperbolic limit cycles.

From these three hyperbolic limit cycles, we can choose two limit cycles, $y_0^{(1)}(b_*)$ and $y_0^{(2)}(b_*)$ such that  $y_0^{(1)}(b_*)<y_0^{(2)}(b_*)$, $d'(y_0^{(1)}(b_*); b_*)<0$, $d'(y_0^{(2)}(b_*); b_*)>0$. By Theorem \ref{dwithb}, for $|b-b_*|\ll 1$ and $b>b_*$,  there exist hyperbolic limit cycles  $y_0^{(1)}(b)$ and $y_0^{(2)}(b)$, where $y_0^{(1)}(b)> y_0^{(1)}(b_*)$
and $y_0^{(2)}(b)< y_0^{(2)}(b_*)$. That is, we can define two curves $y_0=y_0^{(1)}(b)$ and $y_0=y_0^{(2)}(b)$ on $(b_*, b_*+\varepsilon)$, which satisfying that
$d(y_0^{(1)}(b); b)=d(y_0^{(2)}(b); b)\equiv 0$, where  $y_0=y_0^{(1)}(b)$ is increasing and $y_0=y_0^{(2)}(b)$ is decreasing. The above  process can be repeated and we suppose that $\alpha$ is the biggest value such that  $y_0=y_0^{(1)}(b)$ is well defined on $(b_*, \alpha)$ and $\beta$ is the biggest value such that  $y_0=y_0^{(2)}(b)$ is well defined on $(b_*, \beta)$.
Note that $y_0^{(1)}(b)$ and $y_0^{(2)}(b)$ are both hyperbolic limit cycles and $y_0^{(1)}(b)<y_0^{(2)}(b)$.

Obviously, $\{(y_0^{(1)}(b), b) \,|\, b\in (b_*, \alpha)\}, \{(y_0^{(2)}(b), b) \,|\, b\in (b_*, \beta)\}\subset \Omega$. Furthermore, since the two curves are both monotone,
both $y_0^{(1)}(\alpha)\triangleq \lim\limits_{b\rightarrow \alpha-}y_0^{(1)}(b)$  and $y_0^{(2)}(\beta)\triangleq \lim\limits_{b\rightarrow \beta-}y_0^{(2)}(b)$ exist and
$$d(y_0^{(1)}(\alpha); \alpha)=\lim_{b\rightarrow \alpha-}d(y_0^{(1)}(b); b)=0, \quad d(y_0^{(2)}(\beta); \beta)=\lim_{b\rightarrow \beta-}d(y_0^{(2)}(b); b)=0.$$

There are only two possibilities.

(1) $(y_0^{(1)}(\alpha), \alpha)\in \Omega$  or  $(y_0^{(2)}(\beta), \beta)\in \Omega$, that is, $y_0^{(1)}(\alpha)$ or $y_0^{(2)}(\beta)$ is a limit cycle.

If $y_0^{(1)}(\alpha)$ or $y_0^{(2)}(\beta)$ is a hyperbolic limit cycle, then one of the two curves can be well defined in
$(b_*, \alpha+\varepsilon)$ or $(b_*, \beta+\varepsilon)$, where $0<\varepsilon\ll1,$ this is a contraction with the assumption that $\alpha$ and $\beta$ are the biggest value such that
the two curves are well defined. If $y_0^{(1)}(\alpha)$ or $y_0^{(2)}(\beta)$ is a non-hyperbolic limit cycle, then it must be inner stable, this is
a contradiction with the fact that system \eqref{NNE} has only inner unstable non-hyperbolic limit cycle. Thus this possibility cannot occur either.

(2) $(y_0^{(1)}(\alpha), \alpha), (y_0^{(2)}(\beta), \beta)\in \partial\Omega$.

Obviously, we have $(y_0^{(1)}(\alpha), \alpha), (y_0^{(2)}(\beta), \beta)\in \{y_0=b\}$. This contracts with the fact $d(b;b)>0$ for any $b$. Thus this possibility cannot occur, either.
\epf

\begin{Theorem}\label{inunstable} Under the assumptions that $\gamma_L>0>\gamma_R$, $b>0$ and Assumption \ref{parameter-condition} holds. If $\gamma_L+\gamma_R<0$, then system \eqref{NNE} has at most two limit cycles.
\end{Theorem}
\noindent\bpf   The idea of the proof is similar to that of Theorem \ref{instable}, hence we omit some details.

If the conclusion does not hold,  then there exists a  $b_*\in(0,b_M)$ such that system \eqref{NNE} has at least three hyperbolic limit cycles. From these three hyperbolic limit cycles, we can choose two limit cycles,   $y_0^{(1)}(b_*)$ and $y_0^{(2)}(b_*)$ such that  $y_0^{(1)}(b_*)<y_0^{(2)}(b_*)$, $d'(y_0^{(1)}(b_*); b_*)>0$, $d'(y_0^{(2)}(b_*); b_*)<0$.  By Theorem \ref{dwithb}, for $|b-b_*|\ll 1$ and $b<b_*$,  there exist hyperbolic limit cycles  $y_0^{(1)}(b)$ and $y_0^{(2)}(b)$, where $y_0^{(1)}(b)> y_0^{(1)}(b_*)$
and $y_0^{(2)}(b)< y_0^{(2)}(b_*)$. Similarly, we can define two curves $y_0=y_0^{(1)}(b)$ and $y_0=y_0^{(2)}(b)$ on $(\alpha, b_*)$ and $(\beta, b_*)$, which satisfying that
$d(y_0^{(1)}(b); b)=d(y_0^{(2)}(b); b)\equiv 0$, where  $y_0=y_0^{(1)}(b)$ is decreasing and $y_0=y_0^{(2)}(b)$ is increasing. Here we also suppose
that  $\alpha$ and $\beta$ are  the smallest values such that  $y_0=y_0^{(1)}(b)$ and  $y_0^{(2)}(b)$ are well defined on $(\alpha, b_*)$ and $(\beta, b_*)$ respectively.

Obviously, $\{(y_0^{(1)}(b), b) \,|\, b\in (\alpha, b_*)\}, \{(y_0^{(2)}(b), b) \,|\, b\in (\beta, b_*)\}\subset \Omega$. Furthermore, since the two curves are both monotone,
both $y_0^{(1)}(\alpha)\triangleq \lim\limits_{b\rightarrow \alpha+}y_0^{(1)}(b)$  and $y_0^{(2)}(\beta)\triangleq \lim\limits_{b\rightarrow \beta+}y_0^{(2)}(b)$ exist and
$$d(y_0^{(1)}(\alpha); \alpha)=\lim_{b\rightarrow \alpha+}d(y_0^{(1)}(b); b)=0, \quad d(y_0^{(2)}(\beta); \beta)=\lim_{b\rightarrow \beta+}d(y_0^{(2)}(b); b)=0.$$

Similarly, $(y_0^{(1)}(\alpha), \alpha), (y_0^{(2)}(\beta), \beta)\in \partial\Omega$. Obviously, $(y_0^{(1)}(\alpha), \alpha), (y_0^{(2)}(\beta), \beta)\in \{b=0\}\cup \{y_0=\frac{\alpha_R}{\gamma_R-1}+b\}$.  This implies that $d(y_0, b)$ has at least two zeros (counted with multiplicity) on $\{b=0\}\cup \{y_0=\frac{\alpha_R}{\gamma_R-1}+b\}$. This is impossible.
\epf

\vspace{2mm}

For the case $b<0$, of course we assume that $\gamma_L>0>\gamma_R$ and Assumption \ref{parameter-condition} holds. Recall that the function $d(y_0;b)$ is well defined on the interval $(y_0^m, y_0^M)$, where $y_0^m=\max\{0, b\}=0$ and $y_0^M=\min\{\frac{\alpha_R}{\gamma_R-1}+b,\frac{\alpha_L}{\gamma_L+1}\}.$
Since $y_0^m$ must be less than $<y_0^M$ for $d(y_0;b)$ to be well defined, it follows that $b>b_m$.

Now we look $d(y_0;b)$ as a binary function in $(y_0, b)$, which is well defined on the region $\Omega=\{(y_0, b)\,|\, b_m<b<0, 0<y_0<y_0^M\}$. The graph of $\Omega$ has two possibilities.

\begin{figure}[htpb!]
\centering
\subfigure[$\frac{\alpha_L}{\gamma_L+1}\geq\frac{\alpha_R}{\gamma_R-1}$]{
\includegraphics[width=0.45\textwidth]{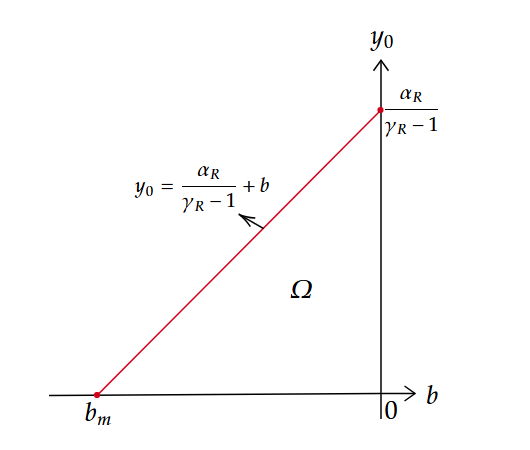}}
\subfigure[$\frac{\alpha_L}{\gamma_L+1}<\frac{\alpha_R}{\gamma_R-1}$]{
\includegraphics[width=0.45\textwidth]{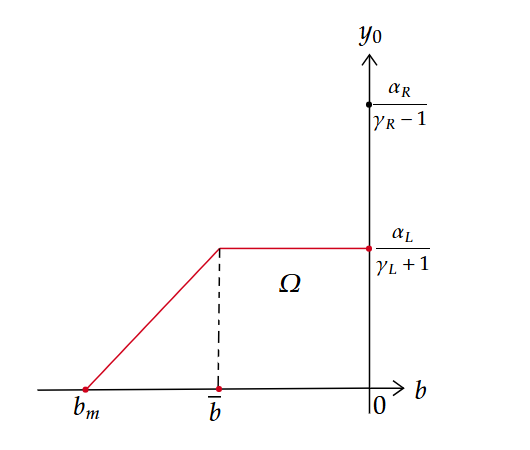}}
\caption{The graph of $\Omega$ for $b<0$}
\label{figy0b<0}
\end{figure}

$(1)$ When $\frac{\alpha_L}{\gamma_L+1}\geq\frac{\alpha_R}{\gamma_R-1},$ we have $y_0^M=\frac{\alpha_R}{\gamma_R-1}+b$ (see Fig. \ref{figy0b<0}$(a)$).

$(2)$ When $\frac{\alpha_L}{\gamma_L+1}<\frac{\alpha_R}{\gamma_R-1},$ there exists a unique $b=\bar{b},$ we have $y_0^M=\frac{\alpha_L}{\gamma_L+1}$ for $b\geq\bar{b}$ and $y_0^M=\frac{\alpha_R}{\gamma_R-1}+b$ for $b<\bar{b}$ (see Fig. \ref{figy0b<0}$(b)$).

Then we will study the sign of $d(y_0; b)$ on some boundaries of $\Omega$. On the boundary $\{b=0\}$, we have known that $d(y_0;0)$ has at most one zero.

On the  boundary $\{y_0=\frac{\alpha_R}{\gamma_R-1}+b\}$, we have
$$\lim\limits_{y_0\to \frac{\alpha_R}{\gamma_R-1}+b}d(y_0;b)=\lim\limits_{y_0\to \frac{\alpha_R}{\gamma_R-1}+b}P_R^{-1}(y_0;b)-P_L\Big(\frac{\alpha_R}{\gamma_R-1}+b\Big)=-\infty<0.$$

\begin{Theorem}\label{rl+rrnot0} Under the assumptions that $\gamma_L>0>\gamma_R$, $b<0$ and Assumption \ref{parameter-condition} holds. If $\gamma_L+\gamma_R\neq0$, then system \eqref{NNE} has at most two limit cycles.
\end{Theorem}

The proof of this theorem is similar to the case $b>0$, which we ignore here for simplicity.

\begin{Theorem}
\label{rl+rr=0}
Under the assumptions that $\gamma_L>0>\gamma_R$ and Assumption \ref{parameter-condition} holds. If $\gamma_L+\gamma_R=0$, then system \eqref{NNE} has at most one limit cycle. If the limit cycle exists, it is stable when $b>0$ and unstable when $b<0$.
\end{Theorem}

Due to the necessity of certain results from the proof of Theorem \ref{dy02}, the proof of this theorem will be presented in Subsection \ref{s4.3}.

\vspace{0.2mm}
Based on the series of theorems provided earlier, Theorem \ref{main th} can be proved. In the following, we present the variation in the number of limit cycles of system \eqref{NNE} with respect to parameter $b$, under a specific parametric condition.

\begin{Theorem}\label{limitcyclewithb}
Under the assumptions that $\gamma_L>0>\gamma_R$, $\Delta_1<0<\Delta_2$ and Assumption \ref{parameter-condition} holds, there exists a $0<\tilde{b}=\tilde{b}(\alpha_{L,R},\gamma_{L,R})<b_M$ such that the following conclusions hold.
\begin{description}
\item $(i)$ When $b\in(-\infty,\bar{b}]\cup(\tilde{b},+\infty),$ system \eqref{NNE} has no limit cycle.
\item $(ii)$ When $b\in(\bar{b},0],$ system \eqref{NNE} has a unique unstable hyperbolic limit cycle.
\item $(iii)$ When $b\in(0,\tilde{b}),$ system \eqref{NNE} has exactly two hyperbolic limit cycles, where the inner limit cycle is stable and the outer one is unstable.
\item $(iv)$ When $b=\tilde{b},$ system \eqref{NNE} has a unique internally stable non-hyperbolic limit cycle.
\end{description}
\end{Theorem}
\noindent\bpf
Given the known conditions, we can deduce that $\gamma_L+\gamma_R<0$, $d''(0;0)<0$ and $d(y_0^M;0)>0$. From Theorem \ref{beq0lc} $iii(b)$, we know that system \eqref{NNE} has a unique unstable limit cycle when $b=0$.

When $b\in(-\infty,b_m)\cup(b_M,+\infty),$ system \eqref{NNE} has no limit cycle obviously since $y_0^m=\min\{0,b\}$ must be less than $y_0^M=\max\{\frac{\alpha_L}{\gamma_L+1}, \frac{\alpha_R}{\gamma_R-1}+b\}$ for $d(y_0;b)$ to be well defined.

We assert that system \eqref{NNE} has at most one limit cycle when $b\in(b_m,0)$. The proof is similar to that of Theorem \ref{instable}. If the assertion does not hold, then there exists a $b^*\in(b_m,0)$ such that system \eqref{NNE} has two hyperbolic limit cycles, as stated in Theorem \ref{main th}. From these two limit cycles $y_0^{(1)}(b^*)$ and $y_0^{(2)}(b^*)$, where $y_0^{(1)}(b^*)<y_0^{(2)}(b^*)$, $d'(y_0^{(1)}(b^*);b^*)>0$ and $d'(y_0^{(2)}(b^*);b^*)<0$ since
$d(b;b)<0$. By Theorem \ref{dwithb}, for $|b-b_*|\ll 1$ and $b<b^*$,  there exist hyperbolic limit cycles  $y_0^{(1)}(b)$ and $y_0^{(2)}(b)$, where $y_0^{(1)}(b)>y_0^{(1)}(b^*)$
and $y_0^{(2)}(b)>y_0^{(2)}(b^*)$. Similarly, we can define two curves $y_0=y_0^{(1)}(b)$ and $y_0=y_0^{(2)}(b)$ on $(\alpha, b^*)$ and $(\beta, b^*)$, which satisfying that
$d(y_0^{(1)}(b); b)=d(y_0^{(2)}(b); b)\equiv 0$, where $y_0=y_0^{(1)}(b)$ is increasing and $y_0=y_0^{(2)}(b)$ is decreasing. Here we also suppose that $\alpha$ and $\beta$ are the smallest values such that  $y_0=y_0^{(1)}(b)$ and  $y_0^{(2)}(b)$ are well defined on $(\alpha,b^*)$ and $(\beta,b^*)$ respectively.

Obviously, $\{(y_0^{(1)}(b), b) \,|\, b\in (\alpha,b^*)\}, \{(y_0^{(2)}(b), b) \,|\, b\in (\beta,b^*)\}\subset \Omega$. Furthermore, since the two curves are both monotone, both $y_0^{(1)}(\alpha)\triangleq \lim\limits_{b\rightarrow \alpha+}y_0^{(1)}(b)$  and $y_0^{(2)}(\beta)\triangleq \lim\limits_{b\rightarrow \beta+}y_0^{(2)}(b)$ exist and
$$d(y_0^{(1)}(\alpha); \alpha)=\lim_{b\rightarrow \alpha+}d(y_0^{(1)}(b); b)=0, \quad d(y_0^{(2)}(\beta); \beta)=\lim_{b\rightarrow \beta+}d(y_0^{(2)}(b); b)=0.$$

Similarly,
$(y_0^{(1)}(\alpha), \alpha), (y_0^{(2)}(\beta), \beta)\in \partial\Omega$. Obviously, $(y_0^{(1)}(\alpha), \alpha), (y_0^{(2)}(\beta), \beta)\in \{y_0=\frac{\alpha_R}{\gamma_R-1}+b\}$.
This implies that $d(y_0;b)$ has at least two zeros (counted with multiplicity) on $\{y_0=\frac{\alpha_R}{\gamma_R-1}+b\}$. This is impossible.

When $b\in(b_m,0)$, we have $y_0^m=0$ and $d(0;b)<d(0;0)=0$ by the proposition $\frac{\partial d(y_0;b)}{\partial b}>0$. And based on previous analysis, we have determined that $y_0^M=\frac{\alpha_R}{\gamma_R-1}+b$ when $b<\bar{b}$, $y_0^M=\frac{\alpha_R}{\gamma_R-1}+\bar{b}=\frac{\alpha_L}{\gamma_L+1}$ when $b=\bar{b}$, and $y_0^M=\frac{\alpha_L}{\gamma_L+1}$ when $b>\bar{b}$. Thus, we have
$$\lim\limits_{y_0\to \frac{\alpha_R}{\gamma_R-1}+b}d(y_0;b)=\lim\limits_{y_0\to \frac{\alpha_R}{\gamma_R-1}+b}P_R^{-1}(y_0;b)-P_L\Big(\frac{\alpha_R}{\gamma_R-1}+b\Big)=-\infty<0$$
when $b\in(b_m,\bar{b})$, by Proposition \ref{propdy0} $(ii)$,
$$\lim\limits_{y_0\to \frac{\alpha_L}{\gamma_L+1}}d(y_0;b)=\lim\limits_{y_0\to \frac{\alpha_L}{\gamma_L+1}}(P(y_0;b)-y_0)=\frac{\alpha_R}{\gamma_R-1}+\bar{b}-\frac{\alpha_L}{\gamma_L+1}=0$$
when $b=\bar{b}$, and
$$\lim\limits_{y_0\to\frac{\alpha_L}{\gamma_L+1}}d(y_0;b)=P_R^{-1}\Big(\frac{\alpha_L}{\gamma_L+1};b\Big)-\lim\limits_{y_0\to \frac{\alpha_L}{\gamma_L+1}}P_L(y_0)=+\infty>0$$
when $b\in(\bar{b},0)$.
According to the fact that system \eqref{NNE} has at most one limit cycle when $b\in(b_m,0)$ and  $d(0;b)<d(0;0)=0$, system \eqref{NNE} has a unique unstable limit cycle for $b\in(\bar{b},0]$ and has no limit cycle for $b\in(b_m,\bar{b})$.

When $b\in(0,b_M)$, by Theorem \ref{pseudohopf}, \ref{dwithb} and \ref{main th}, system
\eqref{NNE} has exactly two limit cycles $y_0=y_0^{(1)}(b)$ and $y_0=y_0^{(2)}(b)$ for a given $0<b\ll1$, where $y_0^{(1)}(b)<y_0^{(2)}(b)$, $d'(y_0^{(1)}(b);b)<0$ and $d'(y_0^{(2)}(b);b)>0$.
According to Theorem \ref{dwithb} $(iii)$, the limit cycles $y_0^{(1)}(b)$ expands and $y_0^{(b)}(b)$ contracts as $b$ increases, then, combine with the fact that system
\eqref{NNE} has no limit cycle for $b\in(b_M,+\infty)$, there must exist a unique $\tilde{b}\in(0,b_M)$ such that at $b=\tilde{b}\in(0,b_M),$ the limit cycles $y_0^{(1)}(b)$ and $y_0^{(2)}(b)$ compose a internally stable semi-stable limit cycle $y_0(\tilde{b})$. Then, system \eqref{NNE} has exactly two limit cycles for $b\in(0,\tilde{b})$ and no limit cycle for $b\in(\tilde{b},b_M)$.

In conclusion, the theorem is proved.
\epf

\section{Proof of Theorems  \texorpdfstring{\ref{dy02}}{} and \texorpdfstring{\ref{rl+rr=0}}{}}
\label{s4}
Based on the expression of $d(y_0;b)$ provided in \eqref{dy0} and the parameter expressions in \eqref{ly0}-\eqref{ry1},  it can be rigorously established that the condition $d(y_0;b)=0$ is equivalent to the following equations hold.
\begin{equation}
\label{y0ty1t}
\left\{
\begin{array}{ll}
\frac{-\alpha_L\psi(v_L,\gamma_L)}{v_L^{\gamma_L-1}(\gamma_L^2-1)(1-v_L^2)}=\frac{-\alpha_R\psi(v_R,\gamma_R)}{v_R^{\gamma_R-1}(\gamma_R^2-1)(1-v_R^2)}+b,\\
\frac{\alpha_Lv_L^{\gamma_L}\psi(v_L^{-1},\gamma_L)}{(\gamma_L^2-1)(1-v_L^2)}=\frac{\alpha_Rv_R^{\gamma_R}\psi(v_R^{-1},\gamma_R)}{(\gamma_R^2-1)(1-v_R^2)}+b,\\
\end{array}
\right.
\end{equation}
where $\alpha_L>0>\alpha_R$, $\gamma_L>1$, $\gamma_R<-1$, $v_L=e^t>1,$ $0<v_R=e^s<1$. The function $\psi(v,\gamma)$ is defined as
$$\psi(v,\gamma)=(\gamma-1)v^{\gamma+1}-(\gamma+1)v^{\gamma-1}+2,$$
under the parameter conditions, it can be easily verified that $\psi(v,\gamma)>0$ holds.

By subtracting the first equation from the second in \eqref{y0ty1t}, we derive the equation
\begin{equation}\label{ML=MR}
  M(\alpha_L,v_L,\gamma_L)=M(\alpha_R,v_R,\gamma_R),
\end{equation}
where the function $M$ is defined as $$M(\alpha,v,\gamma)=\frac{2\alpha(v^{\gamma+1}-1)(v^\gamma-v)}{v^\gamma\big(\gamma^2-1\big)\big(1-v^2\big)}.$$

What's more, according to the parameter expressions in \eqref{ly0}-\eqref{ry1}, we can deduce that
$$d'(y_0;b)=\frac{\mathrm{d}y_1^R}{\mathrm{d}y_0^R}-\frac{\mathrm{d}y_1^L}{\mathrm{d}y_0^L}
=\frac{\mathrm{d}M(\alpha_R,v_R,\gamma_R)}{\mathrm{d}y_0^R}-\frac{\mathrm{d} M(\alpha_L,v_L,\gamma_L)}{\mathrm{d}y_0^L}=F(v_R,\gamma_R)-F(v_L,\gamma_L),$$
where the function $F$ is given by
$$F(v,\gamma)=\frac{\frac{\frac{\partial M}{\partial u}\gamma u}{v}+\frac{\partial M}{\partial v}}{\frac{\frac{\partial y_0}{\partial u}\gamma u}{v}+\frac{\partial y_0}{\partial v}}=\frac{(\gamma-1)(u^2v^2+1)-\gamma(u^2+v^2)+4uv-u^2-v^2}{\gamma  v^2-2uv+v^2-\gamma+1}$$
with $u=v^\gamma$.

If there exists a $y_0^*\in\big(y_0^m,y_0^M\big)$ such that $d(y_0^*;b)=0$ for a given $b,$ we can conclude that the equation \eqref{ML=MR} holds. Then, the second derivative of $d(y_0;b)$ with respect to $y_0$ at $y_0=y_0^*$ can be expressed as
$$d''(y_0^*;b)=\frac{G(v_R,\gamma_R)-G(v_L,\gamma_L)}{M},$$
where the function $G$ is defined as
$$G(v,\gamma)=\frac{\frac{\frac{\partial F}{\partial u}\gamma u}{v}+\frac{\partial F}{\partial v}}{\frac{\frac{\partial y_0}{\partial u}\gamma u}{v}+\frac{\partial y_0}{\partial v}}=\frac{2u(v^2-1)^2(\gamma^2-1)(\gamma uv^2-u^2v-\gamma u+v)(uv-1)(-v+u)}{v(\gamma v^2-2uv+v^2-\gamma+1)^3}$$
with $u=v^\gamma$.

Therefore, to prove Theorem \ref{dy02}, we only need to analyze the sign of $G(v_R,\gamma_R)-G(v_L,\gamma_L)$ when the condition $F(v_R,\gamma_R)=F(v_L,\gamma_L)$ holds.

\subsection{Preliminaries for the proof of Theorem  \texorpdfstring{\ref{dy02}}{}}
In this subsection, we consider $u,$ $v$ and $\gamma$ as independent variables, subject to the constraints $u>v>1$ and $\gamma>1.$

\begin{Lemma}
\label{R1}
The polynomial $R_1(u,v)$ defined by
\begin{align*}
R_1(u,v)=&-u^6 (u^2+1) v^{20}-3 u^5 (5 u^4-2 u^2+5) v^{19}-3 u^4 (u^2+1) (16 u^4-45 u^2\\
&+16)v^{18}-u^3(17 u^8+594 u^6-2526 u^4+594 u^2+17) v^{17}+u^4 (u^2+1) \\
&(53 u^2+14 u-53)(53 u^2-14 u-53) v^{16}+12 u^3 (4 u^8+3773 u^6\\
&-8726 u^4+3773 u^2+4) v^{15}-u^2 (u^2+1)(428 u^8-17291 u^6\\
&+41082 u^4-17291 u^2+428) v^{14}+u(17 u^{12}-14970 u^{10}-472783 u^8\\
&+931608 u^6-472783 u^4-14970 u^2+17) v^{13}-u^2 (u^2+1)\\
&(7713 u^8+948069 u^6-1882520 u^4+948069 u^2+7713) v^{12}\\
&-u (105 u^{12}+397793 u^{10}+1044840 u^8-2824512 u^6\\
&+1044840 u^4+397793 u^2+105) v^{11}+(u^2+1) (17 u^{12}-55874 u^{10}-1803982 u^8\\
&+3675182 u^6-1803982 u^4-55874 u^2+17) v^{10}-u (105 u^{12}+397793 u^{10}\\
&+1044840 u^8-2824512 u^6+1044840 u^4+397793 u^2+105)v^9-u^2 (u^2+1)\\
&(7713 u^8+948069 u^6-1882520 u^4+948069 u^2
+7713)v^8+u (17 u^{12}\\
&-14970 u^{10}-472783 u^8+931608 u^6-472783 u^4-14970 u^2+17) v^7\\
&-u^2 (u^2+1) (428 u^8-17291 u^6+41082 u^4-17291 u^2+428) v^6\\
&+12 u^3 (4 u^8+3773 u^6-8726 u^4+3773 u^2+4) v^5+u^4 (u^2+1)\\
&(53 u^2+14 u-53)(53 u^2-14 u-53) v^4-u^3 (17 u^8+594 u^6-2526 u^4\\
&+594 u^2+17) v^3-3 u^4 (u^2+1) (16 u^4-45u^2+16)v^2-3 u^5 (5 u^4\\
&-2 u^2+5) v-u^6 (u^2+1)
\end{align*}
is strictly negative for $1<v<u\leq u^*$ and has a unique root $v=v(u)\in(1,u)$ for a given $u>u^*,$ where $u^*\approx75.5$ is the unique root of $R_1(u,1)$.
\end{Lemma}
\noindent\bpf
Firstly, we can determine that $R_1(u,1)$ has a unique root $u^*$ by utilizing the Maple built-in command 'sturm'. Then, using Maple built-in command 'fsolve', we can deduce that
\begin{equation*}
R_1(u,1)
\left\{
\begin{array}{lll}
<0,~~if~1<u<u^*,\\
=0,~~if~u=u^*\approx75.5,\\
>0,~~if~u>u^*,
\end{array}
\right.
\end{equation*}
where $R_1(u,1)$ is defined in Appendix.

Furthermore, $R_1(u,u)$ (see Appendix) has no root for any $u>1$ by using the Maple built-in command 'sturm'. It's easy to check that $R_1(u,u)<0$ for any $u>1$.

Using the Maple built-in command 'RealRootIsolate' (see \cite{Dai-Zhao-IJBC2018}), we conclude that the equations $R_1=0$ and $\frac{\partial R_1}{\partial v}=0$ can not hold simultaneously for any $1<v<u.$

$(i)$ To ascertain the sign of $R_1(u,v)$ for any $1<v<u$, we first consider the case $1<u< u^*$. Fixing $u=74<u^*$, we find that $R_1(74,v)$ (see Appendix) has no root for $1<v<74$ by utilizing the Maple built-in command 'sturm'. It's easy to check that $R_1(74,v)<0$ for any $1<v<74$.

\begin{figure}[htpb!]
\centering  
\includegraphics[width=0.4\textwidth]{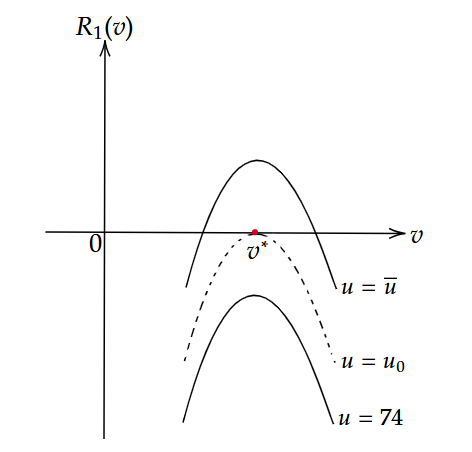}
\caption{Diagram of the function $R_1(u,v)$ for a given $1<u<u^*$.}
\label{fig6}
\end{figure}

We assert that $R_1(u,v)<0$ for all $1<v<u<u^*$. Supposing the contrary, for fixed $\bar{u}\in(1,u^*),$ there exist a $v=v^*\in(1,\bar{u})$, such that $R_1(\bar{u},v^*)\geq0$. There must exist a $u_0\in(1,u^*)$ such that $R_1(u_0,v)$ and $\frac{\partial R_1}{\partial v}(u_0,v)$ have a common root as shown in Fig \ref{fig6} when $u$ varying from $74$ to $\bar{u},$  which leads to a contradiction. Thus, the assertion holds.

$(ii)$ When $u>u^*,$ fixing $u=76>u^*$, we find that $R_1(76,v)$ (see Appendix) has a unique root $v^*$ by utilizing the Maple built-in command 'sturm'. Similar to case $(i)$, this indicates that $R_1(u,v)$ has a unique root $v=v(u)\in(1,u)$ for a given $u>u^*$.

$(iii)$ When $u=u^*$, we can deduce that $R_1(u^*,v)\leq0$ by the proof of case $(i)$. If $R_1(u^*,v)$ has a root, it must be double. This contracts with the fact that $R_1=0$ and $\frac{\partial R_1}{\partial v}=0$ can not hold simultaneously for any $1<v<u.$ Thus, $R_1(u^*,v)<0$ for any $1<v<u^*$.

This completes the proof.
\epf

\begin{Lemma}
\label{R2}
The polynomial $R_2(u,v)$ defined by
\begin{align*}
R_2(u,v)=&-64 u^4 v^{12}+264 u^3 (u^2+1) v^{11}-3 u^2 (688 u^4-1493 u^2+688) v^{10}+2 u (u^2+1)\\
&(68 u^4-1877 u^2+68) v^9+3 u^2 (39521 u^4-78382 u^2+39521) v^8-6 u (u^2+1) (6106 u^4\\
&-10251 u^2+6106) v^7+(1836 u^8-1581246 u^6+3161950 u^4-1581246 u^2+1836) v^6\\
&-6 u (u^2+1)(6106 u^4-10251 u^2+6106) v^5+3 u^2 (39521 u^4-78382 u^2+39521) v^4\\
&+2 u (u^2+1) (68 u^4-1877 u^2+68) v^3-3 u^2 (688 u^4-1493 u^2+688) v^2\\
&+264 u^3 (u^2+1) v-64 u^4
\end{align*}
has a unique root $v=v(u)\in(1,u)$ for a given $u>u^*$, where $u^*$ is defined in Lemma \ref{R1}.
\end{Lemma}
\noindent\bpf
Firstly,  by utilizing the Maple built-in command 'sturm', we can determine that $R_2(u,1)$ (see Appendix) has no root when $u\in[60,+\infty)$. Then, it's easy to check that $R_2(u,1)>0$ for any $u>u^*$.

Furthermore, $R_2(u,u)$ (see Appendix) has no root for any $u>1$ by using the Maple built-in command 'sturm'. Then, it's easy to check that $R_2(u,u)<0$ for any $u>1$, which indicates that $R_2(u,v)$ has at least one root $v=v(u)\in(1,u)$ for a given $u>u^*$.

Using the Maple built-in command 'RealRootIsolate', we conclude that the equations $R_2=0$ and $\frac{\partial R_2}{\partial v}=0$ can not hold simultaneously for any $1<v<u$.

Fixing $u=76>u^*$, we find that $R_2(76,v)$ (see Appendix) has a unique root $v^*$  by
utilizing the Maple built-in command ¡¯sturm¡¯. Thus, by a similar approach to Lemma \ref{R1}, due to the continuous dependence of $R_2(u,v)$ on the parameter $u$, we deduce that $R_2(u,v)$ has a unique simple root $v=v(u)\in(1,u)$ for a given $u>u^*$.
\epf

\begin{Lemma}
\label{H3}
The polynomial $H_3(u,v,\gamma)$ defined by
\begin{align*}
H_3(u,v,\gamma)=&u^3 (v^2-1)^5(2 u^2 v+u v^2+u+2 v)\gamma^5-u^2 v (v^2-1)^4(u^2-1)(u^2 v+14 u v^2+14 u\\
&+v)\gamma^4-2 u (v^2-1)^3(u^6 v^3-8 u^5 v^4-9 u^4 v^5+u^3 v^6-8 u^5 v^2-47 u^4 v^3\\
&+79 u^3 v^4-9 u^2 v^5-9 u^4 v+79 u^3 v^2-47 u^2 v^3-8 u v^4+u^3-9 u^2 v-8 u v^2+v^3)\gamma^3\\
&+2 u v (v^2-1)^2 (u^2-1)(u^4 v^4-15 u^3 v^5-u^2 v^6+u^4 v^2-98 u^3 v^3+127 u^2 v^4\\
&-15 u v^5-15 u^3 v+127 u^2 v^2-98 u v^3+v^4-u^2-15 u v\\
&+v^2)\gamma^2+u (v^2-1)(2 u^6 v^7+16 u^5 v^8-4 u^4 v^9+u^3 v^{10}\\
&+108 u^6 v^5-112 u^5 v^6-34 u^4 v^7-3 u^3 v^8-4 u^2 v^9+2 u^6 v^3-112 u^5 v^4\\
&-164 u^4 v^5+322 u^3 v^6-34 u^2 v^7+16 u v^8+16 u^5 v^2-34 u^4 v^3\\
&+322 u^3 v^4-164 u^2 v^5-112 u v^6+2 v^7-4 u^4 v-3 u^3 v^2-34 u^2 v^3-112 u v^4+108 v^5\\
&+u^3-4 u^2 v+16 u v^2+2 v^3)\gamma-v^2 (u^2-1) (2 u^5 v^7+u^4 v^8+32 u^6 v^4-50 u^5 v^5\\
&+28 u^4 v^6-20 u^3 v^7+u^2 v^8-50 u^5 v^3+38 u^4 v^4-12 u^3 v^5+28 u^2 v^6\\
&+2 u v^7+2 u^5 v+28 u^4 v^2-12 u^3 v^3+38 u^2 v^4-50 u v^5+u^4-20 u^3 v\\
&+28 u^2 v^2-50 u v^3+32 v^4+u^2+2 u v)
\end{align*}
has exact one sign-changing zero point $\gamma=\gamma(u,v)>1$ for any given $(u,v)$, satisfies $1<v<u$ and $u>u^*$.
\end{Lemma}
\noindent\bpf
By the expression of $H_3(u,v,\gamma),$ we have
$$H_3(u,v,1)=-32v(u^3v^5+2u^2v^4+2uv+1)(u-v)^5<0$$
for any $1<v<u$, and
$$H_3(u,v,\gamma)\sim u^3(v^2-1)^5(2 u^2 v+u v^2+u+2v)\gamma^5,~\gamma\to+\infty.$$

Hence, $H_3(u,v,\gamma)$ has at least one zero root $\gamma=\gamma(u,v)>1$. Furthermore, the resultant of $H_3$ and $\frac{\partial H_3}{\partial \gamma}$ with respect to $\gamma$ is
$$Res\big(H_3,\frac{\partial H_3}{\partial\gamma},\gamma\big)=
-4096(v^2-1)^{25}u^{11}v^2(2u^2v+uv^2+u+2v)(uv-1)^{20}
(u-v)^{20}R_2(u,v).$$

By Lemma \ref{R2}, the curve $v=v(u)\in(1,u)$ (satisfying $R_2(u,v(u))=0$ for a given $u>u^*$) divides the $u-v$ plane into two regions:
$$Q_1: \{(u,v)|R_2(u,v)>0,1<v<u,u>u^*\}$$
and
$$Q_2: \{(u,v)|R_2(u,v)<0,1<v<u,u>u^*\}.$$

By evaluating $H_3(u,v,\gamma)$ at points $(u,v)=(76,3)\in Q_1$ and $(u,v)=(76,4)\in Q_2,$ we obtain that $H_3(76,3,\gamma)$ (see Appendix) has a unique root $\gamma\approx36$ and $H_3(76,4,\gamma)$ (see Appendix) has a unique root $\gamma\approx26.$ Therefore, we can conclude that $H_{3}(u,v,\gamma)$ has a unique root $\gamma=\gamma(u,v)>1$ for given $(u,v)$ belongs to regions $Q_1$ and $Q_2$.

For the points on the curve $v=v(u)$, we claim that $H_3(u,v,\gamma)$ has at most one sign-changing zero point $\gamma=\gamma(u,v(u))$ for a given $u>u^*$ and $v(u)\in(1,u)$. If not, we assume that $H_3(u,v,\gamma)$ has two sign-changing zero points $\gamma^{(1)}(u,v(u))$ and $\gamma^{(2)}(u,v(u))$ for a given $u>u^*$ and $v(u)\in(1,u)$. Then,  in a sufficiently small neighborhood of the point $(u,v(u))$, i.e., for $(u,v)\in Q_1$ or $(u,v)\in Q_2$, the function $H_3(u,v,\gamma)$ would have two sign-changing zero points, which contradicts the fact that $H_3(u,v,\gamma)$ has a unique root $\gamma(u,v)$ for given $(u,v)\in Q_1$ or $(u,v)\in Q_2.$

Therefore, $H_{3}(u,v,\gamma)$ has exact one sign-changing zero point $\gamma=\gamma(u,v)>1$ for any given $(u,v)$, satisfies $1<v<u$ and $u>u^*$.

\subsection{The proof of Theorem  \texorpdfstring{\ref{dy02}}{}}
\label{s4.2}
Firstly, we will demonstrate in two steps that the following equations
\begin{equation}
\label{fg}
\left\{
\begin{array}{ll}
F(v_L,\gamma_L)=F(v_R,\gamma_R)\\
G(v_L,\gamma_L)=G(v_R,\gamma_R)
\end{array}
\right.
\end{equation}
do not hold under the specific conditions $v_L,\gamma_L>1,$ $0<v_R<1$, $\gamma_R<-1$ and $\gamma_L+\gamma_R\neq0$.

Step $1$. We claim that $F(v_L,\gamma_L)=F(v_R,\gamma_R)$ determines a unique implicit function $v_L=v_L(\gamma_L,v_R,\gamma_R)>1$ for given $\gamma_L>1$, $0<v_R<1$, $\gamma_R<-1$ and $\gamma_L+\gamma_R\neq0$.

According to the expressions of $F(v_L,\gamma_L)$ and $F(v_R,\gamma_R)$,
$$\lim\limits_{v_L\rightarrow1^+}F(v_L,\gamma_L)-F(v_R,\gamma_R)=-2-F(v_R,\gamma_R)=-\frac{F_1}{F_2},$$
where $F_1=((\gamma_R-1)v_R^2-\gamma_R-1)v_R^{2\gamma_R}+(\gamma_R+1)v_R^2-\gamma_R+1$ and $F_2=\gamma_R v_R^2-2v_R^{\gamma_R+1}+v_R^2-\gamma_R+1$.

The negativity of $F_2$ is affirmed since $\frac{\partial F_2}{\partial v_R}=2(\gamma_R+1)(v_R-v_R^{\gamma_R})>0$ and $F_2(1,\gamma_R)=0$.

Then, focusing on $F_1$, we have $F_1(1,\gamma_R)=\frac{\partial F_1}{\partial v_R}(1,\gamma_R)=\frac{\partial^2 F_1}{\partial v_R^2}(1,\gamma_R)=0$, and
$$\frac{\partial^3 F_1}{\partial v_R^3}=\frac{4\gamma_R(\gamma_R^2-1)v_R^{2\gamma_R}\big(2\gamma_R(v_R^2-1)+v_R^2+1\big)}{v_R^3}<0$$
for any $0<v_R<1, \gamma_R<-1$.

Then, we can conclude that $\frac{\partial^2 F_1}{\partial v_R^2}$ is decreasing with respect to $v_R$. Combined with the fact that $\frac{\partial^2 F_1}{\partial v_R^2}(1,\gamma_R)=0,$ we have $\frac{\partial^2 F_1}{\partial v_R^2}>0$ for any $0<v_R<1, \gamma_R<-1$. Similarly, we can deduce that $\frac{\partial F_1}{\partial v_R}<0$ and $F_1(v_R,\gamma_R)>0$ for any $0<v_R<1, \gamma_R<-1$.

Thus,
$$\lim\limits_{v_L\to1^+}F(v_L,\gamma_L)-F(v_R,\gamma_R)>0$$
for any $\gamma_L>1$, $0<v_R<1$, $\gamma_R<-1$ and $\gamma_L+\gamma_R\neq0$.

In addition,
$$F(v_L,\gamma_L)-F(v_R,\gamma_R)\sim-\frac{\gamma_L-1}{2}v_L^{\gamma_L-1}-F(v_R,\gamma_R),~v_L\to+\infty.$$
Therefore, $F(v_L,\gamma_L)=F(v_R,\gamma_R)$ has at least one root for $v_L(\gamma_L,v_R,\gamma_R)\in(1,+\infty).$

Next, we consider
$$\frac{\partial F(v_L,\gamma_L)}{\partial v_L}=\frac{2v_L^{\gamma_L}(\gamma_L^2-1)(v_L^2-1)F_n}{v_LF_d^2},$$
where $F_d=\gamma_L v_L^2-2v_L^{\gamma_L+1}+v_L^2-\gamma_L+1$
and $F_n=\gamma_L v_L^{\gamma_L+2}-v_L^{2\gamma_L+1}-\gamma_L v_L^{\gamma_L}+v_L.$

We can deduce that
$$F_d(v_L,\gamma_L)<0$$
for any $v_L,\gamma_L>1$ since $\frac{\partial F_d}{\partial v_L}=2(\gamma_L+1)(v_L-v_L^{\gamma_L})<0$
and $F_d(1,\gamma_L)=0.$

In addition, $F_n$ can be rewritten as $-v_Lf,$ where $f=v_L^{2\gamma_L}+\gamma_Lv_L^{\gamma_L-1}-1-\gamma_Lv_L^{\gamma_L+1}$. It can be shown that $f>0$ for $v_L>1$ because $\frac{\partial f}{\partial v_L}=-\gamma_Lv_L^{\gamma_L-2}F_d>0$ and $f(1,\gamma_L)=0.$

Thus, we can deduce that
$$F_n(v_L,\gamma_L)<0$$
for any $v_L,\gamma_L>1$.

Since $\frac{\partial F(v_L,\gamma_L)}{\partial v_L}$ and $F_n$ have the same sign, it follows that $\frac{\partial F(v_L,\gamma_L)}{\partial v_L}<0$ for $v_L>1.$  This completes the proof.

Step $2$. We substitute $v_L(\gamma_L,v_R,\gamma_R)$ that satisfies $F(v_L(\gamma_L,v_R,\gamma_R),\gamma_L)=F(v_R,\gamma_R)$ (here we assume that the values $\gamma_R<-1$ and $0<v_R<1$ are given) into the expression $G(v_L,\gamma_L)-G(v_R,\gamma_R)$. Then, we differentiate this expression with respect to $\gamma_L$, this gives us
$$\frac{\mathrm{d} G}{\mathrm{d} \gamma_L}=\frac{\frac{\partial F}{\partial v}\frac{\partial G}
{\partial\gamma}-\frac{\partial F}{\partial\gamma}\frac{\partial G}{\partial v}}{\frac{\partial F}{\partial v}}\triangleq\frac{H(v,\gamma)}{\frac{\partial F}{\partial v}},$$
where we have omitted the subscript $L$ in the variables for simplicity.

We will prove that the function
$$H(v,\gamma)=\frac{\partial F}{\partial v}\frac{\partial G}
{\partial\gamma}-\frac{\partial F}{\partial\gamma}\frac{\partial G}{\partial v}=H_1\mathrm{ln}v+H_2$$
is sign-definite under the conditions $v,\gamma>1$, where the expressions for $H_1$ and $H_2$ will be given detailed subsequently.

The positivity of $H_2$ is affirmed since
$$H_2(v,\gamma)=-\frac{4u^2(\gamma^2-1)(v^2-1)^2(u-v)(uv-1)}{v^2F_d^5}H_{21}(v,\gamma),$$
where
\begin{align*}
H_{21}(v,\gamma)=&uv(v^2-1)^2(u^2-1)\bigg(\gamma-\frac{(uv^2+u-2v)
(2uv-v^2-1)}{v(v^2-1)(u^2-1)}\bigg)^2\\
&+\frac{(uv^2+u-2v)(2uv-v^2-1)(uv-1)^2(u-v)^2}{v(u^2-1)}>0
\end{align*}
and $u=v^\gamma$.

Therefore, if $H_1=0,$ then we have $H(v,\gamma)=H_2(v,\gamma)>0$ for any $v,\gamma>1.$ If $H_1\neq0$, let's analyze the derivative of the function $\frac{H}{H_1}=\mathrm{ln}v+\frac{H_2}{H_1}$ with respect to $v$. The derivative is expressed as
$$\frac{ \partial \frac{H}{H_1}}{\partial v}=\frac{1}{v}+\frac{\partial\frac{H_2}{H_1}}{\partial v}+\frac{\frac{\partial\frac{H_2}{H_1}}{\partial u}\gamma u}{v}=-\frac{(v^2-1)^2(uv-1)(u-v)F_nH_3(v,\gamma)}{vH_{11}^2},$$
where $H_3(v,\gamma)$ is defined in Lemma \ref{H3} and $u=v^\gamma$.

The Taylor expansion of $H_3(v,\gamma)$ at $v=1$ is given by
$$H_3(v,\gamma)= -\frac{8}{3}\gamma(\gamma^2-1)^5(v-1)^{11}+o\big((v-1)^{11}\big)$$
indicating a leading negative term. Additionally, the asymptotic behavior of $H_3$ as $v\to+\infty$ is analyzed for different ranges of $\gamma$
\begin{equation*}
H_3(v,\gamma)\sim
\left\{
\begin{array}{lll}
-2(\gamma^2-1)(\gamma-1)v^{9+7\gamma},~~if~\gamma<3,v\to+\infty,\\
-64v^{30},~~~~~~~~~~~~~~~~~~~~~~~if~\gamma=3,v\to+\infty,\\
-32v^{6+8\gamma},~~~~~~~~~~~~~~~~~~~~if~\gamma>3,v\to+\infty.
\end{array}
\right.
\end{equation*}

Fixing $\gamma=2$, we find that $H_3(v,2)$ (see Appendix) has no root for any $v>1$ by using the Maple built-in command 'sturm'. Then, it's easy to check that $H_3(v,2)<0$ for any $v>1$.

In order to analyze whether $H_3(v,\gamma)=0$ and $\frac{\partial H_3}{\partial v}(v,\gamma)=0$ can simultaneously hold under the conditions $v,\gamma>1$, we denote $v^\gamma$ as a new variable $u$ (i.e., $v^\gamma=u$) and consider $H_3$ and $\frac{\partial H_3}{\partial v}$ as trinary polynomials in $u$, $v$ and $\gamma$. The existence of intersection points between $H_3$ and $\frac{\partial H_3}{\partial v}$ is equivalent to
their resultant with respect to $\gamma$ yielding the equation
\begin{equation}
\label{H3vH3}
Res\big(H_3,v\frac{\partial H_3}{\partial v},\gamma\big)=
33554432u^{11}v^6(v^2-1)^{25}(uv-1)^{25}(u-v)^{25}(u^2-1)R_1(u,v)=0.
\end{equation}

$(i)$ When $1<v<u\leq u^*$ and $\gamma>1$, by Lemma \ref{R1}, equation \eqref{H3vH3} does not hold when $1<v<u\leq u^*$, that is to say, $H_3=0$ and $\frac{\partial H_3}{\partial v}=0$ can not simultaneously hold when $1<v<u\leq u^*, \gamma>1$. We note that $u=v^\gamma$, combine with the fact that $H_3(v,2)<0$ for any $v>1$, we have $H_3(v,\gamma)<0$ for $1<v<v^\gamma\leq u^*, \gamma>1$.

$(ii)$ When $1<v<u$, $u>u^*$ and $\gamma>2$, by Lemma \ref{H3}, $H_3(u,v,\gamma)$ has exact one sign-changing zero point $\gamma=\gamma(u,v)>1$ for given $1<v<u$ and $u>u^*$. To obtain the sign of $H_3(v,\gamma),$ we only need to show that if $u=v^\gamma,$ then $\gamma<\gamma(u,v).$

The derivatives of $v^\gamma$ with respect to $v$ are
$$\frac{\partial v^\gamma}{\partial v}=\gamma v^{\gamma-1}>0,~~\frac{\partial^2 v^\gamma}{\partial v^2}=\gamma(\gamma-1) v^{\gamma-2}>0,~~\frac{\partial^3 v^\gamma}{\partial v^3}=\gamma(\gamma-1)(\gamma-2) v^{\gamma-3}.$$

When $\gamma>2$, we have
$$
u=v^\gamma>1+\gamma(v-1)+\frac{\gamma(\gamma-1)}{2}(v-1)^2 \Rightarrow
\gamma<\frac{v-3+\sqrt{v^2+8u-6v+1}}{2(v-1)}\triangleq\gamma_2,
$$
while $v^2+8u-6v+1=(v-3)^2+8(u-1)<(u-3)^2+8(u-1)=(u+1)^2$ for $1<v<u$ and $u>u^*$,
we can deduce that $\gamma_2<\frac{u+v-2}{2(v-1)}\triangleq\overline{\gamma_2}.$

When $\gamma=\overline{\gamma_2},$
$$H_3(u,v,\overline{\gamma_2})=-\frac{(u-v)^5}{32}H_{31}(u,v),$$
where
\begin{align*}
H_{31}(u,v)=&9 u^4 v^7 + u^3 (8 u^2 + 3 u + 54) v^6 + u^2 (22 u^3 - 111 u^2 - 306 u + 594) v^5\\
&+ u (16 u^4 - 183 u^3 + 854 u^2 - 936 u + 72) v^4 - u (4 u^4 + 21 u^3 - 180 u^2\\
&+ 356 u - 24) v^3 - u (8 u^4 - 101 u^3 + 358 u^2 - 344 u - 120) v^2 + (-2 u^5\\
&+ 43 u^4 - 338 u^3 + 1186 u^2 - 1848 u + 1024) v - u (u - 2) (u - 4)^2.
\end{align*}

It is verified that $H_{31}(u,1)>0$ (see Appendix) and $H_{31}(u,u)>0$ (see Appendix) for $u>1.$ Using the Maple built-in command 'RealRootIsolate', we confirm that $H_{31}=0$ and $\frac{\partial H_{31}}{\partial v}=0$ can not hold simultaneously for $1<v<u$. In addition, we fix $u=76>u^*$, and find that $H_{31}(76,v)>0$ (see Appendix) for $1<v<76.$ Then, we can deduced that $H_{31}(u,v)>0$ and $H_3(u,v,\overline{\gamma_2})<0$ for $1<v<u,$ $u>u^*$.

Consequently, we can deduce that $\overline{\gamma_2}<\gamma(u,v)$ and $H_3(u,v,\gamma)\leq0$ for any $v^\gamma>u^*$ with $\gamma>2$.

$(iii)$ When $1<v<u, u>u^*$ and $1<\gamma\leq2,$ since
$R_1(v^2,v)$ (see Appendix) has no root for any $v>1$ by using the Maple built-in command 'sturm'. Then, it's easy to check that $R_1(v^2,v)<0$ for any $v>1$. While $R_1\big(u,u^{\frac{1}{2}}\big)=R_1(v^2,v)$, we can deduce that $R_1\big(u,u^{\frac{1}{2}}\big)<0$ and the equation \eqref{H3vH3} does not hold for any $u>u^*$ .

By Lemma \ref{R1}, when $u>u^*,$ $R_1(u,1)>0$ and $R_1(u,v)$ has a unique root $v=v(u)\in (1,u)$, then, we can obtain that $u^{\frac{1}{\gamma}}\geq u^{\frac{1}{2}}>v(u)$ for $1<\gamma\leq2$ by $R_1\big(u,u^{\frac{1}{2}}\big)<0$. Thus, $R_1\big(u,u^{\frac{1}{\gamma}}\big)<0$ for any $u>u^*$ with $1<\gamma\leq2$. And the equation \eqref{H3vH3} does not hold for any $1<v<u, u>u^*$ and $1<\gamma\leq2$ .

Since $H_3(v,2)<0$ for any $v>1$, then we can deduce that $H_3(v,\gamma)<0$ for any $v^\gamma>u^*$ with $1<\gamma\leq2.$

Above all, we obtain that $H_3(v,\gamma)\leq0$ for any $v,\gamma>1$, and when $H_1\neq0$, $\frac{\partial\frac{H}{H_1}}{\partial v}\leqslant0$ for any $v,\gamma>1$ .

In the following, we prove that $H_1<0$ for any $v,\gamma>1$. Focusing on $H_1,$ we have
$$
H_1=\frac{4v^{2\gamma}(\gamma^2-1)(v^2-1)}{v^2F_d^5}H_{11}(v,\gamma),$$
where
\begin{align*}
H_{11}(v,\gamma)=&8v^{6\gamma+5}+(\gamma^2 - 1)v^{5\gamma+8}-\big( (\gamma + 15)(\gamma + 1)v^6 - (\gamma - 1)(\gamma - 15)v^4 + (\gamma^2 - 1)v^2\big)v^{5\gamma}\\
&+ \big(-4\gamma(\gamma - 1)v^9 + 12(\gamma + 1)^2v^7-16(\gamma^2 - 1)v^5 + 12(\gamma - 1)^2v^3 - 4\gamma(\gamma + 1)v\big)v^{4\gamma}\\
&+ \big(\gamma(\gamma^2 - 1)v^{10} - \gamma(5\gamma^2 + 27)v^8 + 10\gamma(\gamma^2 - 1)v^6-10\gamma(\gamma^2 - 1)v^4 + \gamma(5\gamma^2 + 27)v^2 \\
&- \gamma(\gamma^2- 1)\big)v^{3\gamma}+\big(4\gamma(\gamma + 1)v^9 - 12(\gamma - 1)^2v^7 + 16(\gamma^2 - 1)v^5 - 12(\gamma + 1)^2v^3\\
&+ 4\gamma(\gamma - 1)v\big)v^{2\gamma} + \big(-(\gamma^2 - 1)v^8 + (\gamma - 1)(\gamma - 15)v^6 + (\gamma + 15)(\gamma+ 1)v^4 \\
&- (\gamma^2 - 1)v^2\big)v^\gamma - 8v^5.
\end{align*}

The Taylor expansion of $H_{11}(v,\gamma)$ at $v=1$ is
$$H_{11}(v,\gamma)= 4\gamma(\gamma^2-1)^3(v-1)^7+o\big((v-1)^8\big),$$
indicating a leading positive term. Additionally, the asymptotic behavior of $H_{11}$ as $v\to+\infty$ is analyzed for different ranges of $\gamma$
\begin{equation*}
H_{11}(v,\gamma)\sim
\left\{
\begin{array}{lll}
(\gamma^2-1)v^{8+5\gamma},~~if~\gamma<3,v\to+\infty,\\
16v^{23},~~~~~~~~~~~~if~\gamma=3,v\to+\infty,\\
8v^{6\gamma+5},~~~~~~~if~\gamma>3,v\to+\infty.
\end{array}
\right.
\end{equation*}
Thus, we have $H_1<0$ when $v\to1$ and $v\to\infty$.

We claim that $H_1<0$ for any $v\in(1,+\infty).$ If the conclusion does not hold, we consider $H_1$ as a function in $v$, two cases may occur.

$(1)$ If the function $H_1$ exist an even-multiplicity zero $v_e\in(1,+\infty)$ with $H_1(v_e)=0$, without loss of generality, we can assume that $H_1(v)>0$ within the small neighborhood on both sides of $v_e.$ Then, we have $\lim\limits_{v\to v_e^{\pm}}\frac{H}{H_1}(v)=\ln(v_e^{\pm})+\frac{H_2}{H_1}(v_e^{\pm})=+\infty$ within the small neighborhood on both sides of $v_e,$ the monotonicity of $\frac{H}{H_1}$ with respect to $v$ must be opposite. However, when $H_1\neq0,$ we observe the fact that $\frac{\partial \frac{H}{H_1}}{\partial v}\leq0,$ which evidently leads to a contradiction. Hence, this hypothesis is invalid.

$(2)$ If the function $H_1 $ does not exist even-multiplicity zeros, the number of zeros of $H_1$ must be an even number, and we can denote them as $1<v_1<v_2<v_3<\cdots<+\infty$. Given $H_1(v_1)=H_1(v_2)=0$, without loss of generality, we can assume that $H_1(v)>0$ when $v\in(v_1,v_2)$, it follows that $\lim\limits_{v\to v_1}\frac{H}{H_1}(v)=\ln(v_1)+\frac{H_2}{H_1}(v_1)=+\infty$
and $\lim\limits_{v\to v_2}\frac{H}{H_1}(v)= \ln(v_2)+\frac{H_2}{H_1}(v_2)=+\infty.$ Consequently, we can deduce that $\frac{H}{H_1}$ has an even number of even-multiplicity zeros within the interval $v\in(v_1,v_2).$ This contradicts the observed fact that $\frac{\partial \frac{H}{H_1}}{\partial v}\leq0$ for $H_1\neq0$, leading to the invalidation of this hypothesis as well.

Thus, we conclude that $H_1<0$ holds for any $v,\gamma>1.$ Furthermore, since $\lim\limits_{v\to 1^+}\frac{H}{H_1}(v)=0$, we have $H\geq0$ , $\frac{\mathrm{d} G}{\mathrm{d} \gamma}=\frac{H(v,\gamma)}
{\frac{\partial F}{\partial v}}\leq0$ and $\frac{\mathrm{d} G}{\mathrm{d} \gamma}\not\equiv0$ within a small neighborhood of $\gamma$.

Based on the above analysis, if there exists a $y_0^*\in(y_0^m,y_0^M)$, such that $d(y_0^*;b)=0$ and $d'(y_0^*;b)=0$. Then, for a fixed value of $\gamma_R,$ when $\gamma_L=-\gamma_R$, we have $v_L=v_R^{-1}$ by $d(y_0^*;b)=0$ and $d'(y_0^*;b)=0$. Substitute $\gamma_L=-\gamma_R$ and $v_L=v_R^{-1}$ into the expression of $d''(y_0^*;b)$, we have $d''(y_0^*;b)=0$ and $\frac{\mathrm{d} G(v_L,\gamma_L)}{\mathrm{d} \gamma_L}=0$. When $\gamma_L<-\gamma_R$, we can deduce that $d''(y_0^*;b)>0$ since $\frac{\mathrm{d} G(v_L,\gamma_L)}{\mathrm{d} \gamma_L}<0$. Similarly, we have $d''(y_0^*;b)<0$ when $\gamma_L>-\gamma_R$. Above all, we can say that equations \eqref{fg} do not hold under the specific conditions $v_L,\gamma_L>1,$ $0<v_R<1$, $\gamma_R<-1$ and $\gamma_L+\gamma_R\neq0$.

This completes the proof.
\epf

\subsection{The proof of Theorem  \texorpdfstring{\ref{rl+rr=0}}{}}
\label{s4.3}
According to the proof of Theorem \ref{dy02}, if system \eqref{NNE} exists a non-hyperbolic limit cycle, i.e., there exists a $y_0^*\in\big(y_0^m,y_0^M\big),$ such that $d(y_0^*;b)=0$ and $d'(y_0^*;b)=0$ for a given $b.$ When $\gamma_L+\gamma_R=0,$ we can deduce that $d'(y_0^*;b)=0$ if and only if $v_L=v_R^{-1}$. Then we substitute the conditions $\gamma_L+\gamma_R=0$ and $v_L=v_R^{-1}$ into the equation \eqref{ML=MR}, we have $\alpha_L+\alpha_R=0.$ While when $\gamma_L+\gamma_R=0$ and $\alpha_L+\alpha_R=0$, equations \eqref{y0ty1t} does not hold for $b\neq0$.
Thus, system \eqref{NNE} only has hyperbolic limit cycles.

The idea of the proof is similar to that of Theorem \ref{instable}, hence we omit some details. If the conclusion does not hold, then there exists a $b_*\in(0,b_M)$ such that system \eqref{NNE} has two hyperbolic limit cycles. Here we only give the proof for $b>0$, because the proof for $b<0$ is similar, and for the sake of simplicity we will not give its specific proof here.

From these two hyperbolic limit cycles $y_0^{(1)}(b_*)$ and $y_0^{(2)}(b_*)$, where $y_0^{(1)}(b_*)<y_0^{(2)}(b_*)$, we have $d'(y_0^{(1)}(b_*); b_*)<0$ and $d'(y_0^{(2)}(b_*); b_*)>0$ since $d(b;b)>0$. By Theorem \ref{dwithb}, for $|b-b_*|\ll 1$ and $b>b_*$,  there exist hyperbolic limit cycles  $y_0^{(1)}(b)$ and $y_0^{(2)}(b)$, where $y_0^{(1)}(b)>y_0^{(1)}(b_*)$
and $y_0^{(2)}(b)< y_0^{(2)}(b_*)$. Similarly, we can define two curve $y_0=y_0^{(1)}(b)$ and $y_0=y_0^{(2)}(b)$ on $(b_*,\alpha)$ and $(b_*,\beta)$, which satisfying that
$d(y_0^{(1)}(b); b)=d(y_0^{(2)}(b); b)\equiv 0$, where  $y_0=y_0^{(1)}(b)$ is increasing and $y_0=y_0^{(2)}(b)$ is decreasing. Here we also suppose that  $\alpha$ and $\beta$ are the largest values such that  $y_0=y_0^{(1)}(b)$ and  $y_0^{(2)}(b)$ are well defined on $(b_*,\alpha)$ and $(b_*,\beta)$ respectively.

Obviously, $\{(y_0^{(1)}(b), b) \,|\, b\in (b_*,\alpha)\}, \{(y_0^{(2)}(b), b) \,|\, b\in (b_*,\beta)\}\subset \Omega$. Furthermore, since the two curves are both monotone,
both $y_0^{(1)}(\alpha)\triangleq \lim\limits_{b\rightarrow \alpha-}y_0^{(1)}(b)$  and $y_0^{(2)}(\beta)\triangleq \lim\limits_{b\rightarrow \beta-}y_0^{(2)}(b)$ exist and
$$d(y_0^{(1)}(\alpha); \alpha)=\lim_{b\rightarrow \alpha-}d(y_0^{(1)}(b); b)=0, \quad d(y_0^{(2)}(\beta); \beta)=\lim_{b\rightarrow \beta-}d(y_0^{(2)}(b); b)=0.$$

Similarly,
we have $(y_0^{(1)}(\alpha), \alpha), (y_0^{(2)}(\beta), \beta)\in \partial\Omega$. Obviously, we have $(y_0^{(1)}(\alpha), \alpha), (y_0^{(2)}(\beta), \beta)\in \{y_0=b\}$.
This contracts with the fact $d(b; b)>0$ for any $b$.

\section{Discussion}
\label{s5}
In this paper, we have proven that system \eqref{NNE} has at most two limit cycles, and the multiplicity of any limit cycle of such a system does not exceed two. This conclusion was reached through a detailed analysis of the roots of the successor function $d(y_0;b)$, along with an examination of the monotonic variation of $d(y_0;b)$ with respect to the parameter $b$ and the corresponding variation in limit cycles. Notably, the proof methodology employed in this study has the potential to be extended to systems of the form \eqref{GE} that exhibit $SS$, $SN,$ $SIN,$ and $NIN$ types, which we intend to explore in our future research. Moreover, it is worth highlighting that Li et al., in their work \cite{Li-Chen-QTDS2020}, have proven the absence of sliding limit cycles in system \eqref{NNE}. This observation indicates that even when considering the possibility of sliding limit cycles, the total number of limit cycles of system \eqref{NNE} remains unchanged.

In addition, when proving the multiplicity of limit cycles in this paper, we utilized the method of symbolic computation, which involves a relatively complicated process. Exploring how to use a simpler method to prove this aspect is also worth studying.

\section*{Acknowledgements}

This work is supported by the National Natural Science Foundation of China, under Grant Number 12171491.

\section*{Appendix}

\begin{align*}
R_1(u,v=1)=&(u+1)^2\big(u^2+u+1\big)\big(17u^{10}-227u^9-71526u^8-610029u^7-1615317u^6\\
&+4554960u^5-1615317u^4-610029u^3-71526u^2-227u+17\big).
\end{align*}

\begin{align*}
R_1(u,v=u)=&-81u^6\big(u^2+1\big)\big(u^{20}-25u^{18}-456u^{16}+25326u^{14}-7797u^{12}\\
&-31194u^{10}-7797u^8+25326u^6-456u^4-25u^2+1\big).
\end{align*}

\begin{align*}
R_1(u=74,v)=&-899358946693952v^{20}-998033287229652576v^{19}\\
&-236273294846159095872v^{18}-6233874231058218081704v^{17}\\
&+13828786452622076515408v^{16}+20501425245268386910464v^{15}\\
&-11457482312591678745030912v^{14}+28434278502998776616331514v^{13}\\
&-212676989914124169685024388v^{12}-354523933163399054340696386v^{11}\\
&+994839798467750981296989933v^{10}-354523933163399054340696386v^9\\
&-212676989914124169685024388v^8+28434278502998776616331514v^7\\
&-11457482312591678745030912v^6+20501425245268386910464v^5\\
&+13828786452622076515408v^4-6233874231058218081704v^3\\
&-236273294846159095872v^2-998033287229652576v-899358946693952.
\end{align*}

\begin{align*}
R_1(u=76,v)=&-1113227487383552v^{20}-1268771824564122624v^{19}\\
&-308489817962516201472v^{18}-8356334350772692924096v^{17}\\
&+18055403762437022103808v^{16}+27280979403158247822336v^{15}\\
&-15784652469256215920728512v^{14}+40621920223278737678705036v^{13}\\
&-292552429520779710016454288v^{12}-490771670006273317238659564v^{11}\\
&+1560081942219797860180360833v^{10}-490771670006273317238659564v^9\\
&-292552429520779710016454288v^8+40621920223278737678705036v^7\\
&-15784652469256215920728512v^6+27280979403158247822336v^5\\
&+18055403762437022103808v^4-8356334350772692924096v^3\\
&-308489817962516201472v^2-1268771824564122624v\\
&-1113227487383552.
\end{align*}

\begin{align*}
R_1(u,v=2)=&17408u^{14}-127360u^{13}-97804288u^{12}-1143554440u^{11}-5615140720u^{10}\\
&-5209018030u^9+5416299723u^8+11881429500u^7+5416299723u^6\\
&-5209018030u^5-5615140720u^4-1143554440u^3-97804288u^2-127360u+17408.
\end{align*}

\begin{align*}
R_2(u,v=1)=&1836u^8-73000u^7-1348248u^6+43032u^5+2700488u^4+43032u^3\\
&-1348248u^2-73000u+1836.
\end{align*}

$$R_2(u,v=u)=-216u^4(8u^{12}-393u^{10}+8490u^8-15968u^6+8490u^4-393u^2+8).$$

\begin{align*}
R_2(u=76,v)=&-2135179264v^{12}+669494588928v^{11}-397583235316224v^{10}\\
&+1982571342746336v^9+22839237334338480v^8-536478278900935632v^7\\
&+1738931358905378380v^6-536478278900935632v^5+22839237334338480v^4\\
&+1982571342746336v^3-397583235316224v^2+669494588928v-2135179264.
\end{align*}

\begin{align*}
H_3(u=76,v=3,\gamma)=&509522997149696\gamma^5-11464849789747200\gamma^4\\
&-254059618640269312\gamma^3+47622937841740800\gamma^2\\
&+3611925187779284480\gamma-24809993679363631200.
\end{align*}

\begin{align*}
H_3(u=76,v=4,\gamma)=&15836668279200000\gamma^5-278265426516000000\gamma^4\\
&-3292245295420608000\gamma^3-216103625559360000\gamma^2\\
&+33077473098545767680\gamma-141369532357567852800.
\end{align*}

\begin{align*}
H_3(v,\gamma=2)=&-32v^{22}+\big(-6v^9+302v^7-218v^5+18v^3\big)v^{14}+\big(23v^{10}-1020v^8+1226v^6\\
&+4v^4-297v^2\big)v^{12}+\big(-32v^{11}+1658v^9-2802v^7-202v^5+1778v^3\\
&-432v\big)v^{10}+\big(18v^{12}-1352v^{10}+3290v^8-3290v^4+1352v^2-18\big)v^8\\
&+\big(432v^{11}-1778v^9+202v^7+2802v^5-1658v^3+32v\big)v^6+\big(297v^{10}-4v^8\\
&-1226v^6+1020v^4-23v^2\big)v^4+\big(-18v^9+218v^7-302v^5+6v^3\big)v^2+32v^6.
\end{align*}

$$H_{31}(u,v=1)=32u^5-160u^4+96u^3+800u^2-1600u+1024.$$

\begin{align*}
H_{31}(u,v=u)=&u\big(17u^{10}+25u^9-41u^8-493u^7+1419u^6-657u^5-599u^4+29u^3+1316u^2\\
&-1880u+1056\big).
\end{align*}

\begin{align*}
H_{31}(u=76,v)=&300259584v^7+20407994240v^6+51947461024v^5+34832612448v^4\\
&-10765745952v^3-17069780576v^2-3778140160v-29154816.\\
\end{align*}

\begin{align*}
R_1(u=v^2,v)=&-9v^{10}\big(51v^{26}-94v^{25}+682v^{24}+471v^{23}+2699v^{22}+33297v^{21}\\
&+38096v^{20}-13039v^{19}+197314v^{18}-252196v^{17}+309782v^{16}-354726v^{15}\\
&+192534v^{14}-296674v^{13}+192534v^{12}-354726v^{11}+309782v^{10}-252196v^9\\
&+197314v^8-13039v^7+38096v^6+33297v^5+2699v^4+471v^3+682v^2\\
&-94v+51\big)(v+1)^2.
\end{align*}

\end{document}